\documentclass[ss,noinfoline]{imsart}

\RequirePackage[OT1]{fontenc}
\RequirePackage{amsthm,amsmath}
\RequirePackage[numbers]{natbib}
\RequirePackage[colorlinks,citecolor=blue,urlcolor=blue]{hyperref}
\usepackage{amsfonts}
\usepackage{bbm}
\usepackage{graphicx}
\usepackage{tikz-cd}
\usepackage{pict2e}

\startlocaldefs
\numberwithin{equation}{section}
\theoremstyle{plain}

\theoremstyle{definition}
\newtheorem{example}{Example}[section]
\newtheorem*{remark}{Remark}
\newtheorem*{remarks}{Remarks}
\endlocaldefs

\begin{document}

\begin{frontmatter}

%
%
%

{\huge{\bf Note:} {\it Corrected and extended parts of this paper will soon be published as a book of the same title}}
\title{\!\!\!\!\!\!\!\!\!\!\!\!Causality from the Point of View of Statistics\!\!\!\!\!\!\!\!\!\!\!\!\!}
\runtitle{Causality: Viewpoint of Statistics}

\begin{aug}
\author{\fnms{Jos\'{e} A.} \snm{Ferreira}\ead[label=e1]{jose.ferreira@rivm.nl}}

\address{Department of Statistics, Informatics and Modelling\\
National Institute for Public Health and the Environment (RIVM)\\ 
Antonie van Leeuwenhoeklaan 9,
3721 MA Bilthoven,
The Netherlands\\
\printead{e1}}

\runauthor{Author}


\end{aug}

\begin{abstract}
We present a basis for studying questions of cause and effect in statistics 
which subsumes and reconciles the models proposed by Pearl, Robins, Rubin and others, and which, as far as  mathematical notions and notation are concerned, is entirely conventional. In particular, we show that, contrary to what several authors had thought, standard probability can be used
to treat problems that involve notions of causality, and in a way not essentially different from the way it has been used in the area generally known (since the 1960s, at least) as `applied probability'. 
Conventional, elementary proofs are given of some of the most important results obtained by the various schools of `statistical causality', and a variety of examples considered by those schools are worked out in detail.
Pearl's `calculus of intervention' is examined anew, and its first two rules are formulated and proved by means of elementary probability for the
first time since they were stated 25 years or so ago. 
\vspace{-0.2cm}
\end{abstract}

%
%
%
%
%


\begin{keyword}
\kwd{Causality}
\kwd{Markovian models}
\kwd{Identifiability}
\kwd{Confounding\!\!}
\vspace{-0.1cm}
\end{keyword}
\vspace{-0.2cm}
\end{frontmatter}

\newpage
\section{Background}\label{Background}
Let $\Omega=[0,1]$, ${\cal F}$ be the family of Borel subsets of $\Omega$, $\mathbf{P}$ Lebesgue measure on $\Omega$,
and consider the probability space $(\Omega,{\cal F},\mathbf{P})$. In this space one can compute probabilities of events such as $\{\omega\}$,
\[
E:=\bigcup_{n}[a_n,b_b],\quad \Omega\cap \mathbb{Q}\equiv \{\omega\in\Omega: \omega\in\mathbb{Q}\},
\]
or
\[
N_j:=\left\{\omega\equiv\sum_{m\ge 1}\frac{\omega_m}{10^m}\in\Omega: \lim_{n\longrightarrow\infty}\frac{1}{n}\sum_{m=1}^n\delta_{\omega_m,j}=\frac{1}{10}\right\}
\]
($j=0,1,...,9$). Indeed, $\mathbf{P}(\{\omega\})=0$, so that the probability of drawing any particular element of $\Omega$ is 0, 
$\mathbf{P}(E)=\sum_{n} \mathbf{P}\left([a_n,b_b]\right)$ provided the subintervals of $\Omega$ in the union are disjoint,
$\mathbf{P}(\Omega\cap \mathbb{Q})=0$, so that the event `$\omega$ is rational' has probability 0, and $\mathbf{P}(N_j)=1$ for all $j$---Borel's \emph{normal number} theorem.
We interpret each $\omega\in\Omega$ as a possible `random draw' from the unit interval and say that the probability that
the draw falls between $a$ and $b$ is $\mathbf{P}([a,b])=b-a$ ($a\le b$), that it is impossible to draw any particular $\omega$,
that the draw is certain to be irrational, etc.

If $X:\Omega\rightarrow\mathbb{R}$ is a random variable, the real number $X(\omega)$ is called the `realization' of $X$ associated with the random draw $\omega$;
if $A$ is an arbitrary Borel subset of $\mathbb{R}$, the event that a realization of $X$ falls in $A$, defined as $X^{-1}(A):=\{\omega\in\Omega: X(\omega)\in A\}$ and often abbreviated
as $\{X\in A\}$ or `$X$ is in $A$', has probability $\mathbf{P}(X\in A)\equiv \mathbf{P}(X^{-1}(A))$.
In particular, if $U$ is the identity mapping on $\Omega$ then $\mathbf{P}(U\le u)=u$ ($u\in [0,1]$),
$\mathbf{P}(U\mbox{ is irrational})\equiv \mathbf{P}(U\in\Omega\backslash\mathbb{Q})=1$, and so on; $U$ is said to be a \emph{standard uniform} random variable, or a random variable with the standard uniform distribution (the distribution being the mapping $A\rightarrow \mathbf{P}(A)$).

Moreover, if $D_n(\omega)$ denotes the $n$-th digit in the non-terminating decimal expansion of $U(\omega)$, so that
\begin{equation}\label{ExpansionOfUniform}
U(\omega)=\sum_{n\ge 1}\frac{D_n(\omega)}{10^n}\equiv 0\!\cdot\!D_1(\omega)D_2(\omega)D_3(\omega)\cdots,
\end{equation}
then $D_1,D_2,\ldots$ are random variables on $\Omega$ and, for instance,  the probability that the second and fourth digits equal 6 and 3, $\mathbf{P}\left(D_2=6,D_4=3\right)$, is
\[
\sum_{j,k=0}^9\left(\frac{j}{10}+\frac{6}{10^2}+\frac{k}{10^3}+\frac{3}{10^4}< U\le \frac{j}{10}+\frac{6}{10^2}+\frac{k}{10^3}+\frac{4}{10^4}\right)=\frac{1}{10^2},
\]
which, because $\mathbf{P}(D_n=j)=1/10$ for all $j$,
shows that $\mathbf{P}\left(D_2=6,D_4=3\right)=\mathbf{P}(D_2=6)\mathbf{P}(D_4=3)$; and since the digits 3 and 6 play no special role in this calculation it follows that $D_2$ and $D_4$ are independent. More generally, $D_1,D_2,\ldots$ are independent random variables, each {\emph{uniform on}} $\{0,1,\ldots,9\}$.

Thus a standard uniform random variable may be identified with a decimal expansion like (\ref{ExpansionOfUniform}) consisting of independent random digits uniform on $\{0,1,\ldots,9\}$.
But the digits of $U$ can be arranged in a two-dimensional array by the `diagonal method' in such a way that for each $\omega$ the correspondence 
\[
\begin{array}{ccccccccc}
U_1(\omega)\,\, \leftrightarrow\!\!\! &  0\!\cdot\!\!D_1(\omega)& \!\!\!\!\!D_3(\omega) & \!\!\!\!\!D_6(\omega) & \!\!\!\!\!D_{10}(\omega) & \!\!\!\!\!D_{15}(\omega) & \!\!\!\!\!D_{21}(\omega) & \!\!\!\!\!D_{28}(\omega) & \cdots\\
U_2(\omega)\,\, \leftrightarrow\!\!\! &  0\!\cdot\!\!D_2(\omega) & \!\!\!\!\!D_5(\omega) & \!\!\!\!\!D_9(\omega) & \!\!\!\!\!D_{14}(\omega) & \!\!\!\!\!D_{20}(\omega) & \!\!\!\!\!D_{27}(\omega)  & \cdots & \cdots\\
U_3(\omega)\,\, \leftrightarrow\!\!\! & 0\!\cdot\!\!D_4(\omega) & \!\!\!\!\!D_8(\omega) & \!\!\!\!\!D_{13}(\omega) & \!\!\!\!\!D_{19}(\omega) & \!\!\!\!\!D_{26}(\omega) & \cdots &\cdots &\cdots\\
U_4(\omega)\,\, \leftrightarrow\!\!\! & 0\!\cdot\!\!D_7(\omega) & \!\!\!\!\!D_{12}(\omega) & \!\!\!\!\!D_{18}(\omega) & \!\!\!\!\!D_{25}(\omega) & \cdots & \cdots & \cdots & \cdots\\
U_5(\omega)\,\, \leftrightarrow\!\!\! &  0\!\cdot\!\!D_{11}(\omega) & \!\!\!\!\!D_{17}(\omega) & \!\!\!\!\!D_{24}(\omega) & \cdots & \cdots & \cdots & \cdots & \cdots\\
U_6(\omega)\,\, \leftrightarrow\!\!\! &  0\!\cdot\!\!D_{16}(\omega) & \!\!\!\!\!D_{23}(\omega) & \cdots & \cdots & \cdots & \cdots & \cdots & \cdots\\
U_7(\omega)\,\, \leftrightarrow\!\!\! &  0\!\cdot\!\!D_{22}(\omega) & \cdots & \cdots & \cdots & \cdots & \cdots & \cdots & \cdots\\
\multicolumn{9}{c}{\dotfill}
\end{array}
\]
determines a sequence $U_1(\omega),U_2(\omega),\ldots$ of numbers in $\Omega$, and thereby (the random digits $D_n$ making up the decimal expansions being independent) a sequence $U_1,U_2,\ldots$ of standard uniform random variables, all defined on the same probability space $(\Omega,{\cal F},\mathbf{P})$. This construction is concrete to the extent that 
numbers such as $\pi-3$ and $e-2$, say, are concrete and that to each of them there corresponds a countable sequence of irrational numbers
$U_1(\omega),U_2(\omega),\ldots$ Whether a particular $U(\omega)=\omega$ is a normal number or not is quite irrelevant (and probably unascertainable);
what is relevant is that a draw from the probability space yields a countable number of independent draws with the same distribution, that sets of draws
possessing certain properties have probability 1 or 0, etc. 

If $F$ is an arbitrary right-continuous distribution function on $\mathbb{R}$ and its `inverse' $F^{-1}$ is defined by $F^{-1}(u)=\min\{x\in\mathbb{R}: F(x)\ge u\}$ then
$\mathbf{P}(F^{-1}(U)\le x)=F(x)$, i.e.~the random variable $X=F^{-1}(U)$ has distribution function $F$. It follows that the correspondence between each $U(\omega)=\omega$ and 
$X_1(\omega),X_2(\omega),\ldots$, where $X_n(\omega)=F^{-1}(U_n(\omega))$ for each $\omega$ and $n$, determines a sequence $X_1,X_2,\ldots$ of independent 
 random variables on the same probability space, each of which has distribution function $F$. More generally, a sequence $X_1,X_2,\ldots$ of independent random variables with
distribution functions $F_1,F_2,\ldots$ is defined by the transformations $X_n=F_n^{-1}(U_n)$. Still more generally, if $F_1$ is a distribution function and, for each $n\in\mathbb{N}$ and each $\mathbf{x}\in\mathbb{R}^{n}\!$, $F_{n+1,\mathbf{x}}\!$ is a distribution function in $\mathbb{R}$ then $X_1\!=\!F_1^{-1}(U_1)$ and
\[
X_{n+1}=F_{n+1,\mathbf{X}_n}^{-1}(U_{n+1}),\quad \mathbf{X}_n=(X_1,\ldots,X_{n-1},X_n),\quad n\in\mathbb{N},
\]
defines a sequence of \emph{generally dependent} random variables on $(\Omega,{\cal F},\mathbf{P})$ with joint distribution functions given by $G_1=F_1$,
\[
\begin{array}{ccc}
G_{n}(\mathbf{x}_n)&\!\!:=\!\!\!& \!\!\!\!\!\!\!\!\mathbf{P}(X_1\le x_1,X_2\le x_2,\ldots,X_n\le x_n)\\
&\!\!=\!\!\!\!&\int_{-\infty}^{x_{n-1}}\cdots \int_{-\infty}^{x_1} F_{n,\mathbf{x'_{n-1}}}(x_{n})\mbox{d}G_{n-1}(\mathbf{x'_{n-1}}),
\end{array}
\]
where $\mathbf{x}_n=(x_1,x_2,\ldots,x_n)\in\mathbb{R}^n$, $n=2,3,\ldots$ If for instance $F_{n+1,\mathbf{x}}$ does not depend on $n$ and if it depends on $\mathbf{x}\in\mathbb{R}^{n}$ only through the last coordinate then the sequence $X_1,X_2,\ldots$ is a homogeneous Markov chain.
Evidently, the mapping of each draw $U(\omega)=\omega$ to each sequence $X_1(\omega),X_2(\omega),\ldots$ is not just a simple
rearrangement of the digits of $\omega$, as is the case with $U_1(\omega),U_2(\omega),\ldots$; but, again, the unfolding of a number like $\pi-3$ diagonally
into a two-dimensional array followed by the transformation of the rows of the array by a real-valued function is hardly less concrete than $\pi-3$ itself.

Finally, one can define an uncountable set $Y=(Y_t)_{t\ge 0}$ of
random variables---a random function or `continuous-time' stochastic process on $\mathbb{R}_+$---on the same probability space. 
Indeed, the procedure used above to get $U_1,U_2,\ldots$ from $U$ can be applied to each of the $U_n$s to yield an independent sequence of sequences $(V_{n,m})_{m\ge 1}$ of independent standard uniform random variables, and these in turn can be transformed into an independent sequence of sequences $(X_{n,m})_{m\ge 1}$ of independent random variables
with any given distribution functions. But then, for example, two such sequences, say $(X_{1,m})_{m\ge 1}$ and $(X_{2,m})_{m\ge 1}$, the first of which may be assumed strictly positive, suffice to determine for each $\omega$ a right-continuous step function\footnote{As usual, ${\bf 1}_E(\omega)$ equals 1 for $\omega\in E$ and 0 for $\omega\not\in E$.}
$Y_t(\omega)=\sum_{m\ge 1}X_{2,m}(\omega){\bf 1}_{\{X_{1,1}+\cdots+X_{1,m}\le t\}}(\omega)$,
whose jumps occur at the points
$X_{1,1}(\omega)+\cdots+X_{1,m}(\omega)$ and have sizes $X_{2,m}(\omega)-X_{2,m-1}(\omega)$;
in particular, if the first sequence has the same exponential distribution and the second is degenerate at 1 (i.e.~\!$\mathbf{P}(X_{2,n}=1)=1$ for all $n$) then $Y$ is a Poisson process. 

For a second example based on all the independent sequences $(X_{n,m})_{m\ge 1}$ obtained above, consider for each $n\in\mathbb{N}$ some real-valued functions $(f_{m,n})_{m\ge 1}$ on $\mathbb{R}_+$ such that $f_{m,n}(t)=0$ for $t\not\in \,\,]n-1,n]$ and set
\[
Y^{(n)}_t(\omega)=\sum_{m\ge1}f_{m,n}(t) X_{n,m}(\omega)
\]
for $t\!\in\![n-1,n]$ and $\omega\in\Omega$ for which the series converges and $Y^{(n)}_t(\omega)\!=\!0$ for other $(t,\omega)$.  
Under certain conditions on the $f_{m,n}$s and on the $X_{n,m}$s (e.g.~\!\cite{ItoNisio1968}),
for each $n$ the function $t\rightarrow Y^{(n)}_t(\omega)$ is continuous on $[n-1,n]$ for $\omega$ in a set of probability 1.
Under such conditions, for $\omega$ in a set of probability 1 the function 
\[
Y_t(\omega)=\sum_{n\ge1}\left(Y^{(n)}_t(\omega)+\sum_{j=1}^{n-1}Y^{(j)}_j(\omega)\right)\!{\bf 1}_{[n-1,n[}(t)
\]
is continuous on $\mathbb{R}_+$.
For instance, if the $X_{n,m}$s are standard normal random variables
and $f_{n,m}(t)=\frac{\surd{2}}{\pi m}\left[1-\cos(m\pi\{t-(n-1)\})\right]{\bf 1}_{[n-1,n[}(t)$ then the
$Y=(Y_t)_{t\ge 0}$ thus defined is a Wiener (or Brownian motion) process.

\section{Models of causality}\label{DeterministicModel}

When one says that something affects, is the cause of, or has an effect on something else, then one implies that it does so after a certain moment and in a certain situation which
may help bring about the effect and which, up to that moment, may itself affect and be affected by the cause.
Mathematically, this idea can be expressed by the equations
\begin{equation}\label{EquationDeterministicModel}
r=\rho(\mathbf{x},t),\quad \tau(\mathbf{x}, t)=0,
\end{equation}
where $\mathbf{x}$ describes the situation prior to the moment when the effect comes into force,
$ t$ stands for the cause---called \emph{treatment} here on account of problems of cause and effect typically encountered in applications---$r$ stands for the joint result of, or the \emph{response} to,
$\mathbf{x}$ and $ t$, and $\tau$ and $\rho$ are some functions, $\tau$ being vector-valued in order that the second equation may consist of one or more equations relating $ t$ to components (or coordinates)
of $\mathbf{x}$.\footnote{In principle $\mathbf{x}$, $ t$ and $r$  may be thought of as arbitrary sets of numerical, textual and pictorial elements, but we shall think of them mainly as vectors of real numbers or functions.}

Although the second equation may define $ t$ implicitly as a function of $\mathbf{x}$---and the components of the latter implicitly as functions of its other components
and of $ t$---it really is a \emph{constraint} representing the interaction between $ t$ and $\mathbf{x}$ prior to the moment when the effect comes into force. Thus, for a given $\mathbf{x}$
there may be a single $ t$ satisfying $\tau(\mathbf{x}, t)=0$, there may be none, but at least for one
$\mathbf{x}$ there must be $t_1\neq t_2$ such that $\tau(\mathbf{x}, t_1)=\tau(\mathbf{x}, t_2)=0$---for if $ t$ were uniquely determined by $\mathbf{x}$ then the response would
be a function of the situation alone and there would be no treatment to talk of, a possibility that we exclude.
On the other hand, except in artificial cases neither $\mathbf{x}$ nor $ t$ will be functions of $r$; and even if $(\mathbf{x}, t)$ can be written as a function of $r$ there should be no reason for one to interpret (\ref{EquationDeterministicModel}) as saying that the response causes the situation or the treatment.\footnote{For instance, if 
$\rho(\mathbf{x}, t)=(\mathbf{x}-1)(t-1)+\mathbf{x}(\mathbf{x}-1)/2+t(t+1)/2$ then $\rho$ is one-to-one in $\mathbb{N}\times \mathbb{N}$ and maps this set onto $\mathbb{N}$, so one might very well consider that $r$ causes $(\mathbf{x}, t)$ rather than that $\mathbf{x}$ and $ t$ cause $r$; but of course that would require us to interpret $\mathbf{x}$, $ t$ and $r$ and their relation with ``the moment when the effect comes into force'' in a different way.
The possibility of one confusing cause and effect in our two equations, having read the definitions that presuppose them, is as `problematic' as the possibility of one concluding that the motion of a marble rolling down an inclined plane is the cause of gravity or of the inclination of the plane.}

In order to turn these ideas into a `causal model' we only need to specify a few things about
$\rho$ and $\tau$, such as their domains and codomains, sets $\cal{ X}$ and ${\cal T}$ of situations and treatments, and a set ${\cal A}\subset\cal{ X}\times {\cal T}$ of {\it admissible elements} $(\mathbf{x},t)$ satisfying $\tau(\mathbf{x}, t)=0$ and which can be used to compute values of $r=\rho(\mathbf{x}, t)$.
For, irrespectively of how apt such model is to represent a real-life situation, we may consider trying a definite range of values of $\mathbf{x}$ and $ t$  in (\ref{EquationDeterministicModel}) 
and ask \textit{whether and to what extent the treatment affects the response}---whether the {\it response function} $(\mathbf{x}, t)\rightarrow\rho(\mathbf{x}, t)$ varies with $ t$ for some $\mathbf{x}$, and, if the treatment is shown to have an effect, to try and quantify that effect in terms of differences between responses $\rho(\mathbf{x},t_1)$ and $\rho(\mathbf{x},t_2)$  to different treatments $t_1$, $t_2$ under one or more situations $\mathbf{x}$.

Of course, in applications the response function is to be regarded as unknown (there would be no problem to study otherwise), even if one may assume partial information about it, such as that
$\rho(\mathbf{x}, t)$ is constant in $\mathbf{x}'$ when $\mathbf{x}=(\mathbf{x}',\mathbf{x}'')$
over a range of values of $\mathbf{x}''$; and something about $\tau$ must be known, such as the elements of ${\cal A}$. Thus, in applications of the model one would not know the exact expression of $\rho$, and perhaps not that of $\tau$, but one would see the values of $\rho$ assumed on ${\cal A}$.

And the model applies at least to some concrete problems.
For example, the lifetime of a set of batteries used in a radio receiver depends on the brand of the batteries
and possibly on such factors as elapsed shelf life, temperature and rate of usage during operation, etc.
In principle, by varying the brand of the batteries 
and using the radio in such a way as to keep the other factors fixed, one should be able to determine which brand lasts longer.
Such `studies' or `experiments' can take place in everyday life and need not require much planning or care in order to provide satisfactory answers.
Thus, if two sets of batteries being tested have very different lifetimes then slight differences between the temperature, rate of usage, elapsed shelf life, etc.,
experienced during the two periods of testing will probably not mask the superiority of one set over the other: in the language of (\ref{EquationDeterministicModel}),
if $t_1$ and $t_2$ represent the sets of batteries and $\mathbf{x}_1$ and $\mathbf{x}_2$ the situations during the two periods of testing,
one expects $r_1:=\rho(\mathbf{x}_1,t_1)\approx \rho(\mathbf{x},t_1)$ and $r_2:=\rho(\mathbf{x}_2,t_2)\approx \rho(\mathbf{x},t_2)$
provided $\mathbf{x}_1\approx \mathbf{x}\approx \mathbf{x}_2$, in which case the comparison between $r_1$ and $r_2$ is a suitable replacement
for the more accurate comparison between $\rho(\mathbf{x},t_1)$ and $\rho(\mathbf{x},t_2)$.\footnote{In this example, the allowed combinations of situations and treatments may be represented by a table with binary entries indicating whether they apply or not, and $\tau$ may consist of a vector of sums of products of indicator functions or Kronecker $\delta$s of the entries of the table.}
But the same sort of consideration serves to show that the effect of the factor of interest may be \emph{confounded} with (or by) the effects of incidental factors: if $\rho(\mathbf{x},t_1)\approx \rho(\mathbf{x},t_2)$ for each $\mathbf{x}$ and if $\rho(\mathbf{x}_1, t)$ and $\rho(\mathbf{x}_2, t)$ are very different for $\mathbf{x}_1\neq \mathbf{x}_2$ and each $ t$, then
the large difference between $r_1$ and $r_2$ will be a poor substitute for the small differences between
$\rho(\mathbf{x}_1,t_1)$ and $\rho(\mathbf{x}_1,t_2)$ and between $\rho(\mathbf{x}_2,t_1)$ and $\rho(\mathbf{x}_2,t_2)$.

Let us say in connection with our model that \emph{confounding} may occur or exist, or that the situation may be a \emph{confounder} of the treatment,
whenever $\tau$ is not constant in $\mathbf{x}$ and there exist
$(\mathbf{x}_1,t_1),(\mathbf{x}_2,t_2)\in {\cal A}$
such that $\mathbf{x}_1\neq \mathbf{x}_2$, $t_1\neq t_2$ and
\[
\rho(\mathbf{x}_1,t_i)\neq\rho(\mathbf{x}_2,t_i)\quad\mbox{ for }i=1,2.
\]
While admissibility implies that the choices of $t$ and $\mathbf{x}$ are  subordinated to each other to some extent, the last condition implies that an investigation of the effect of the treatment by comparing 
$\rho(\mathbf{x},t_1)$ with $\rho(\mathbf{x},t_2)$ for a given $\mathbf{x}$ is not
 guaranteed---since on the one hand $(\mathbf{x},t_1)$ and $(\mathbf{x},t_2)$ need not be both admissible, and on the other hand $\rho(\mathbf{x}_1,t_1)$ and $\rho(\mathbf{x}_2,t_2)$ may differ because 
$\mathbf{x}_1\neq \mathbf{x}_2$ rather than because $t_1\neq t_2$.

However, even when confounding is possible one may still be able to study the effect of the treatment 
if there exist at least two different points $(\mathbf{x},t_1), (\mathbf{x},t_2)\in {\cal A}$ with which to compare $\rho(\mathbf{x},t_1)$ and $\rho(\mathbf{x},t_2)$, or if for distinct
$\mathbf{x}_1$ and $\mathbf{x}_2$ such that $(\mathbf{x}_1,t_1), (\mathbf{x}_2,t_2)\in {\cal A}$ one has
$\rho(\mathbf{x}_1,t)=\rho(\mathbf{x}_2,t)$ for each $ t$, which permits comparing 
$\rho(\mathbf{x}_1,t_1)=\rho(\mathbf{x}_2,t_1)$ with $\rho(\mathbf{x}_1,t_2)=\rho(\mathbf{x}_2,t_2)$.
Accordingly, in such cases we shall say that \emph{confounding can be removed}---or that one can \emph{correct for confounding}, or something similar.
In the first case, removal of confounding is achieved by \emph{matching} on the situation, i.e. by taking points $(\mathbf{x},t_1)$ and $(\mathbf{x},t_2)$ in ${\cal A}$ with $t_1\neq t_2$  for each one of as many as possible $\mathbf{x}$, and comparing $\rho(\mathbf{x}, t_1)$ with $\rho(\mathbf{x},t_2)$.
In the second case, using the fact that $\rho(\mathbf{x}, t)=\rho(\mathbf{y},t)$ for certain  
$\mathbf{x}\neq \mathbf{y}$ and $ t$ in a certain range,  
one corrects for confounding by \emph{stratifying} or partitioning a subset of $\cal{ X}$ into the sets
\[
{\cal X}_{\mathbf{x}}=\big\{\mathbf{y}\in{\cal X}:\,\tau(\mathbf{y},t)=\tau(\mathbf{x},t)=0,\,
                                            \rho(\mathbf{y},t)=\rho(\mathbf{x}, t)\mbox{ for some }t\in{\cal T}\big\}
\]
and the complement of their union, and comparing the responses to various $t\in {\cal T}$ within each set of the resulting partition.

It goes without saying that in order to do this in a real-life problem one has to know enough about the response and be allowed to pick favourable combinations of situations and treatments; one may assume wrongly that $\rho(\mathbf{x},t)$ is constant in $\mathbf{x}'$ for $\mathbf{x}=(\mathbf{x}',\mathbf{x}'')$ in a certain set and make the `confounded comparison' between
$\rho(\mathbf{x}',\mathbf{x}'',t_1)$ and $\rho(\mathbf{x}',\mathbf{x}''',t_2)$, and then, in effect,
{\it confounding occurs}.

\begin{example}\label{PedanticExample}
Take some functions
$F,G:\mathbb{R}^2\rightarrow\mathbb{R}$, $H:\mathbb{R}^3\rightarrow\mathbb{R}$,
some intervals $A,B\subset\mathbb{R}$, and put
\[
\,\,\tau(\mathbf{x},t)=F(\mathbf{x})-|t|{\bf 1}_{A}(x_1)-t{\bf 1}_{A^{c}}(x_1),
\quad\quad\quad\quad\quad\quad\quad\quad\quad\quad\quad
\]
\[
\rho(\mathbf{x},t)=G(x_2,t){\bf 1}_{B}(x_2)+H(\mathbf{x},t){\bf 1}_{B^c}(x_2),\quad\,
\mathbf{x}=(x_1,x_2)\in\mathbb{R}^2, t\in\mathbb{R}.
\]

If $x_1\in A^c$ then $\tau(\mathbf{x},t)=0$ $\Leftrightarrow$ $F(\mathbf{x})=t$ and 
$\rho(\mathbf{x},t)=\rho(\mathbf{x},F(\mathbf{x}))$, so in general the last argument of $\rho$ cannot be varied in order to study its effect (though it can be in trivial cases, as when $\rho(\mathbf{x},t)$ depends on $t$ alone and one can take $\mathbf{x}\neq \mathbf{x}'$ such that
$F(\mathbf{x})=t\neq t'=F(\mathbf{x}')$).

If $x_1\in A$ and $x_2\in B^{c}$ then $\rho(\mathbf{x},t)=H(\mathbf{x},t)$ and $\tau(\mathbf{x},t)=0$
$\Leftrightarrow$ $F(\mathbf{x})=|t|$, so there may be $\mathbf{x}$ with which to compare
$H(\mathbf{x},-t)$ and $H(\mathbf{x},t)$; that is, we may study the treatment effect by matching on
values of $\mathbf{x}\in A\times B^{c}$.

If $x_1\in A$ and $x_2\in B$ then $\rho(\mathbf{x},t)=G(x_2,t)$ and
$\tau(\mathbf{x},t)=0$ $\Leftrightarrow$ $F(\mathbf{x})=|t|$, so
there may be $(x_1,x_2)$, $(x_1',x_2)$ satisfying $t=F(x_1,x_2)$ and $-t=F(x_1',x_2)$, with
which we can compare
$G(x_2,-t)$ and $G(x_2,t)$; that is, we may study the treatment effect within strata
${\cal X}_{x_1,x_2}:=\left\{(x_1',x_2): x_1'\in A, F(x_1',x_2)=-F(x_1,x_2)\right\}$,
$(x_1,x_2)\in A\times B$.$\hfill\square$
\end{example}

The possibilities of studying the effect of treatment on the response
depend on ${\cal A}$ and on $\rho$.
If $\tau$ is constant in $ t$ then $\mathbf{x}$ and $ t$ can be chosen at will, i.e.~within the limits imposed by the constraints on the $\mathbf{x}$s and $ t$s but independently of each other.
In our example of the radio batteries, the model with such $\tau$ corresponds to the possibility of carrying out an \emph{experimental study} involving, for instance,
buying batteries of various brands fresh from the factories at the appropriate times and testing each brand by playing the radio continuously until the batteries die,
at constant volume, in a room with controlled temperature, humidity, etc.

If on the contrary $\tau$ varies with $ t$ and $\mathbf{x}$ then the extent to which confounding can be removed depends on how rich
${\cal A}$ is.\footnote{For example, let $n\in\mathbb{N}$ and consider $\tau(x,t)=x^2+t^2-n^2$.
If ${\cal A}\subset\mathbb{Z}^2$ then for each situation $x$ there are at most two treatments that can be used to evaluate the response.
If ${\cal A}\subset\mathbb{R}^2$, two treatments, namely $t_{\pm}=\pm\sqrt{n^2-x^2}$, can be considered for each situation $x\!\in\, ]-n,n[$, and therefore confounding may be removed by comparing
$\rho(x,-t)$ with $\rho(x,t)$ and $\rho(-x,-t)$ with $\rho(-x,t)$;
but even in this case the range of admissible $(x,t)$s is quite restricted compared to $\mathbb{R}^2$.}
In the example of the batteries, this last version of the model may represent an \emph{observational study} in which a careless experimenter buys batteries from various shops at unplanned times, depending on availability and other circumstances, and does not take particular care in fixing  temperature, rate of usage, etc.
If the experimenter carried out many such experiments he might be fortunate enough to have a few of them in which batteries of different brands happened to have the same elapsed shelf life and the rate of usage, temperature, etc., just happened to be constant or exerted no influence on the life of the batteries. 

\section{The simplest statistical model}\label{StatisticalModel}

The simplest statistical model is
\begin{equation}\label{CausalModelII}
R_n=\rho\big(V_n,\mathbf{X}_n,T_n\big),\quad T_n=\tau\big(U_n,\mathbf{X}_n\big),\quad n\in\mathbb{N},
\end{equation}
where $(\mathbf{X}_n)_{n\ge 1}$ is a sequence of random `situations' taking values in ${\cal X}$ (e.g.~in $\mathbb{R}^d$), $(U_n)_{n\ge 1}$ and $(V_n)_{n\ge 1}$ are sequences of random variables
such that $U_n$ and $V_n$ are both standard uniform conditionally on $\mathbf{X}_n$,
all defined on our probability space $(\Omega,{\cal F},\mathbf{P})$,
$(T_n)_{n\ge 1}$ is a sequence of random treatments
taking values in ${\cal T}$ (e.g.~\!in $\mathbb{N}_0$), $(R_n)_{n\ge 1}$ is the corresponding sequence of responses (e.g.~taking values in $\mathbb{R}$),
and $\rho$ and $\tau$ are real-valued functions defined on $[0,1]\times {\cal X}\times{\cal T}$ and $[0,1]\times{\cal X}$, respectively.
Just like in (\ref{EquationDeterministicModel}), there is no danger of confusing cause and effect in equations (\ref{CausalModelII}):
$R_n$ is a function of $(\mathbf{X}_n,T_n)$ and represents an event occurring after the determination of the latter; more precisely, the calculation of $R_n(\omega)$
for each $\omega\in\Omega$ is really preceded by that of $U_n(\omega)$, $V_n(\omega)$ and
$\mathbf{X}_n(\omega)$ followed by that of $T_n(\omega)$ from $\mathbf{X}_n(\omega)$---and if there must be an `ultimate cause of things'
it is the drawing of $\omega$.\footnote{\label{FootnoteConcern}A concern about standard mathematical notation not being able to represent relationships of cause and effect appears to be quite common;
see for instance the text leading to footnote 5 of \cite{Pearl2009a}, the second paragraph on p.~\!291 of \cite{CoxWermuth04}.}

We shall sometimes refer to this as the {\it basic model}, to distinguish it from the more detailed models considered in sections \ref{LevelsOfCausality} and \ref{IdentificationEstimation}.

The admissible elements in the models of section \ref{DeterministicModel} have here a counterpart in the sequence $(\mathbf{X}_n,T_n)_{n\ge 1}$ of pairs of situations and treatments, which possesses more structure and may be endowed with such features as independence or dependence and
stationarity or nonstationarity.

An interesting variant of (\ref{CausalModelII}) is obtained by letting $T_n$ be constant with probability 1, or letting it take
a fixed number of values at fixed values of $n$, while allowing the $\mathbf{X}_n$s to be nondegenerate.
If we regard our model as representing a self-contained `world' or `system' evolving from a random draw $\omega$, we can think of this option as the external forcing of a fixed treatment, or of a fixed pattern of treatments, upon the system---which without
that external \emph{intervention} would have evolved randomly according to its own laws.
Thus, we will see that when a real-life problem is represented by a model
like (\ref{CausalModelII}) it is sometimes useful to consider the associated \emph{intervention model}
\begin{equation}\label{CausalModelIIa}
R_n=\rho\big(V_n,\mathbf{X}_n,t_n\big),\quad n\in\mathbb{N},
\end{equation}
where $(\mathbf{X}_n)_{n\ge 1}$, $(V_n)_{n\ge 1}$ have the same distributions as in (\ref{CausalModelII})
and $(t_n)_{n\ge 1}$ is a sequence of numbers in ${\cal T}$---which in particular may be set equal to a constant in order to represent the forcing of the same treatment
in every situation in which the treatment may exert an effect.\footnote{\label{FootnoteToInterventionModel}When the $t_n$s are all equal, 
the intervention model corresponds to an application of Pearl's `\emph{do} operator'---``which simulates physical interventions by deleting certain functions
[in this case $\tau$] from the model, replacing them by a
constant [...], while keeping the rest of the model unchanged''---to (\ref{CausalModelII}); see p. 107 of \cite{Pearl2009b} or pp. 54-57 of \cite{PearlEtAl2016}.
Of course, in many non-experimental real-life situations the enforcement of a treatment upon a system does affect other aspects of that system (for instance, forcing an
individual to stop smoking may cause him to increase his consumption of alcohol or foodstuffs, or to recalcitrate in unanticipated ways; cf.~\!section 8.2 of
\cite{CoxWermuth04}), but, as we shall see later on, an intervention model, even if purely hypothetical or unrealizable, may, in principle, serve a good purpose.
Note that (\ref{CausalModelIIa}) really defines a \emph{family of models}, one model for each sequence $(t_n)_{n\ge 1}$, and that (\ref{CausalModelIIa}) could also
be seen as a special case of (\ref{CausalModelII}) where $\tau$ depends only on the index $n$ of the pair $(U_n,\mathbf{X}_n)$.}

These models bring with them a distinction between \emph{observable} and \emph{unobservable} random variables: in (\ref{CausalModelII}),
the situations, the treatments and the responses are observable because
they are regarded as empirical data which can be observed or measured in connection with a cause and effect relationship; 
the $U_n$s and $V_n$s, on the other hand, are regarded as unobservable, because they were brought in precisely as `factors' lying beyond
the situations and treatments, i.e.~\!beyond what can be observed or measured. Of course, unobservable random variables are as mathematical and as `real'
as observable random variables;
the main reason for distinguishing
them is that only the observable ones are used in statistical procedures---such as versions of the matching and stratification procedures mentioned in section \ref{DeterministicModel}, or procedures to estimate functions such as
$(\mathbf{x},t)\rightarrow \mathbf{E}[R_n|\mathbf{X}_n=\mathbf{x},T_n=t]$, which plays an  important role in the sequel.

As in section \ref{DeterministicModel}, to say that treatment has an effect on the response---or that there exists a \emph{treatment effect}---in connection with (\ref{CausalModelII}) or (\ref{CausalModelIIa})
is essentially to say that the response function $\rho$ is not constant in its third argument.
However, in order to avoid exceptional cases which are irrelevant for applications, we shall say that treatment has an effect on the response if, for some $n$, conditionally on
${\mathbf X}_n$, $\rho(V_n,{\mathbf X}_n,t)$ varies with $t$ with positive probability; that is, there is a treatment effect
if and only if\footnote{It is customary and sometimes convenient to
refer to the probability distribution $A\rightarrow \mathbf{P}(X\in A)$ of a random variable $X$ as \emph{the law of $X$}, which we denote by ${\cal L}(X)$;
similarly, the law of $Y$ conditional on the event $\{X=x\}$, written ${\cal L}(Y|X=x)$, stands for $A\rightarrow \mathbf{P}(Y\in A|X=x)$.}
\[
{\cal L}\big(\rho(V_n,{\mathbf X}_n,t)\big|{\mathbf X}_n={\mathbf x}\big)\equiv
{\cal L}\big(\rho(V_n,\mathbf{x},t)\big|{\mathbf X}_n={\mathbf x}\big)
\]
\[
\quad\quad\quad\quad\quad\quad\quad={\cal L}\big(\rho(V_n,\mathbf{x},t)\big)
\]
is not constant in $t$ for some $n$ and some ${\mathbf x}$ such that
${\mathbf P}({\mathbf X}_n\in N)>0$ for every neighbourhood $N$ of ${\mathbf x}$.

\begin{remark}\label{FootnoteConditioningCanBeRemoved}
Note that since $U_n$ and $V_n$ are standard uniform conditionally on $\mathbf{X}_n$, each of them is standard uniform unconditionally as well.
It follows that the conditioning on $\{\mathbf{X}_n\!=\!\mathbf{x}\}$ can be removed from
${\cal L}\big(\rho(V_n,\mathbf{x},t)\big|\mathbf{X}_n\!=\!\mathbf{x}\big)$
and from other expressions involving the law of $\rho(V_n,\mathbf{x},t)$, as done
in the last inequality.\,$\square$
\end{remark}

In real-life problems the $(\mathbf{X}_n,T_n,R_n)$s are meant to represent observations on `individuals' (`patients', `units', etc.) and serve as models for the responses of those individuals to treatments under certain conditions. The hypothesis that the treatment has no effect
on the response then means that the response of the $n$-th individual is fully determined by the situation and by incidental factors represented by $V_n$.
On the other hand, it is enough that a single individual's response be affected by the treatment in a
`realizable' situation for the treatment to have an effect.

It follows from all this that a study of the treatment effect on the response amounts to a study of the conditional---or unconditional---laws of the random variables
$\rho(V_n,\mathbf{X}_n,t)$, $t\in{\cal T}$, the
so-called \emph{potential outcomes}, and of related random variables such as
$\rho(V_n,\mathbf{X}_n,t)-\rho(V_n,\mathbf{X}_n,t')$ or
$\rho(V_n,\mathbf{X}_n,t)/\rho(V_n,\mathbf{X}_n,t')$, $t\neq t'$.
But these are unobservable (only $\rho(V_n,\mathbf{X}_n,T_n)$ is observable), so it is not obvious that their laws can be estimated from observed data.
In fact, we shall now see that \emph{the possibility of estimating ${\cal L}\big(\rho(V_n,\mathbf{X}_n,t)\big|\mathbf{X}_n=\mathbf{x}\big)$ from data}---the possibility of studying the effect of treatment on
response---corresponds to \emph{the possibility of removing confounding} in (\ref{StatisticalModel});
for the $\mathbf{X}_n$s are confounders of treatment because unless one is somehow able to `fix' them it is generally impossible to know whether differences in response are due to differences in the situations or to differences in the treatment.

Ideally, the study of the treatment effect in our statistical model (\ref{CausalModelII}) would involve the comparison between the 
values of $\rho\big(V_n,\mathbf{X}_n,t\big)$ for various $t\in{\cal T}$---the potential outcomes of the $n$-th individual.
However, at each draw $\omega$ we only get to see $\rho\big(V_n(\omega),\mathbf{X}_n(\omega),t\big)$
for a single $t$, namely $t=T_n(\omega)$;
we do not see any of its \emph{counterfactuals}, i.e.~the set of $\rho\big(V_n(\omega),\mathbf{X}_n(\omega),t\big)$ for $t\neq T_n(\omega)$---if we did,
then the differences
\[
\rho\big(V_n(\omega),\mathbf{X}_n(\omega),T_n(\omega)\big)-\rho\big(V_n(\omega),\mathbf{X}_n(\omega),t\big)
\]
for $t\neq T_n(\omega)$ would reveal the effect of the treatment and solve our problem.
Even in the intervention model (\ref{CausalModelIIa}) the most that one is given to see and is allowed to compare
(e.g.~\!when the $\mathbf{X}_n$s are discrete) 
at a single draw $\omega$ are pairs of responses
$\rho\big(V_m(\omega),\mathbf{X}_m(\omega),t_m\big)$ and 
$\rho\big(V_n(\omega),\mathbf{X}_n(\omega),t_n\big)$
such that $\mathbf{X}_m(\omega)=\mathbf{X}_n(\omega)=\mathbf{x}$ for some $\mathbf{x}\in{\cal X}$ and $t_m\neq t_n$, which, because they generally differ
in the first argument of $\rho$, will, on their own, seldom tell us whether differences between them are caused by the treatment or by $V_m$ and $V_n$ (whose realizations are unobservable).
It seems evident, then, that in order to study the effect of treatment based on a single realization $(\mathbf{X}_n(\omega),T_n(\omega),R_n(\omega))_{n\ge 1}$
of (\ref{CausalModelII}) one has to make do with the comparison of empirical conditional 
averages such as
\[
\frac{\sum\limits_{n=1}^N \mathbf{1}_{\{\mathbf{X}_n=\mathbf{x},T_n=t\}}(\omega)\rho\big(V_n(\omega),\mathbf{x},t\big)}
{\sum\limits_{n=1}^N \mathbf{1}_{\{\mathbf{X}_n=\mathbf{x},T_n=t\}}(\omega)}
\,\,\mbox{ and  }\,\,\,
\frac{\sum\limits_{n=1}^N \mathbf{1}_{\{\mathbf{X}_n=\mathbf{x},T_n=t'\}}(\omega)\rho\big(V_n(\omega),\mathbf{x},t'\big)}
{\sum\limits_{n=1}^N \mathbf{1}_{\{\mathbf{X}_n=\mathbf{x},T_n=t'\}}(\omega)}
\]
for $t\neq t'$ and various $\mathbf{x}$---and hope that the $V_n(\omega)$s for which $T_n(\omega)=t$ are not essentially different,
as far as their role in $\rho$ is concerned, from those for which $T_n(\omega)=t'$. However, there is nothing in the conditions introduced so far that
prevents us from having, for instance, $V_n(\omega)<1/2$ for $n$ such that $\mathbf{1}_{\{\mathbf{X}_n=\mathbf{x},T_n=t\}}(\omega)=1$ and
$V_n(\omega)\ge 1/2$ for $n$ such that $\mathbf{1}_{\{\mathbf{X}_n=\mathbf{x},T_n=t'\}}(\omega)=1$,
an occurrence that would \emph{confound} the workings of $t$ in $\rho\big(V_n(\omega),\mathbf{x},t\big)$ and could produce a difference between the two averages even
if $\rho$ were constant in its third argument. It is therefore necessary to require, 
for the purpose of studying
the treatment effect in model (\ref{CausalModelII}), that 
$V_n$ be independent of $T_n$ conditionally on $\mathbf{X}_n$, or, what is equivalent, that \emph{$U_n$ and $V_n$ be independent conditionally on
$\mathbf{X}_n$}.
This condition, which we shall refer to as \emph{unconfoundedness}, should make it possible to remove confounding
because each $V_n$ contributing to one average
has the same distribution as each $V_n$ contributing to the other, so that differences between averages that seem too extreme can be correctly attributed to the fact that
$t\neq t'$ rather than to randomness.

When and to what extent it is possible to study the effect of the treatment by comparing
empirical conditional averages also depends on the properties
of those averages, which in turn depend on properties of the sequence $(\mathbf{X}_n,T_n,R_n)_{n\ge 1}$ other than unconfoundedness.
Thus, although we shall not consider questions of testing for and estimating a treatment effect (treated in \cite{Rosenbaum2002}, \cite{Freedman:2010},
\cite{ImbensRubin:2015} and \cite{Ferreira:2015}, for example), we must note that when the distribution of $(\mathbf{X}_n,T_n,R_n)$ is independent of $n$,
and hence is the same as that of a generic vector $(\mathbf{X},T,R)$ defined on the same probability space, empirical conditional averages of the form
\begin{equation}\label{GeneralSampleConditionalAverage}
\frac{\sum_{n=1}^N \mathbf{1}_{\{\mathbf{X}_n=\mathbf{x},T_n=t\}}f(R_n)}
{\sum_{n=1}^N \mathbf{1}_{\{\mathbf{X}_n=\mathbf{x},T_n=t\}}},
\end{equation}
where $f$ is some real-valued function, converge under certain conditions and
in a certain sense to
\begin{equation}\label{GeneralConditionalAverage}
\mathbf{E}[f(R)|\mathbf{X}=\mathbf{x},T=t]=\mathbf{E}[f(R_n)|\mathbf{X}_n=\mathbf{x},T_n=t]
\end{equation}
as $N\rightarrow\infty$ (e.g.~\!with probability 1 if the $(\mathbf{X}_n,T_n,R_n)$s are independent).\footnote{This is generally true if the 
$\mathbf{X}_n$s and $T_n$s are discrete, which we assume for simplicity in most of the paper; if the $\mathbf{X}_n$s are not discrete, the indicators
$\mathbf{1}_{\{\mathbf{X}_n=\mathbf{x},T_n=t\}}$ in (\ref{GeneralSampleConditionalAverage}) are replaced by
$\mathbf{1}_{\{\mathbf{X}_n\in B_{N}(\mathbf{x}),T_n=t\}}$, where $B_{N}(\mathbf{x})$ is a neighbourhood shrinking to $\mathbf{x}$ as $N\rightarrow\infty$.}
In the sequel we shall often refer to the possibility of estimating expectations such as these \emph{from observed data},
by which we imply that the right-hand side of (\ref{GeneralConditionalAverage}) 
is independent of $n$ and can be estimated consistently by (\ref{GeneralSampleConditionalAverage}) with samples---not necessarily
random samples---$(\mathbf{X}_1,T_1,R_1),...,(\mathbf{X}_N,T_N,R_N)$.

In order to confirm that under unconfoundedness the comparison of empirical conditional averages---those of (\ref{GeneralSampleConditionalAverage})
for varying $t$---provides a way of studying the effect of
treatment, consider (\ref{GeneralConditionalAverage}) with $f=\mathbf{1}_{A}$ and varying $A$, i.e.~\!the law
\[
{\cal L}\big(R_n\big|\mathbf{X}_n=\mathbf{x},T_n=t\big)=
{\cal L}\big(\rho(V_n,\mathbf{x},t)\big|\mathbf{X}_n=\mathbf{x},T_n=t\big).
\]
As we have indicated, when based on a single realization of $(\mathbf{X}_n,T_n,R_n)_{n\ge 1}$ the study of the treatment effect 
amounts to the comparison of the conditional laws
\begin{equation}\label{ComparisonOfLawsI}
{\cal L}\big(\rho(V_n,\mathbf{x},t)\big|\mathbf{X}_n=\mathbf{x},T_n=t\big)\quad\mbox{and}\quad
{\cal L}\big(\rho(V_n,\mathbf{x},t')\big|\mathbf{X}_n=\mathbf{x},T_n=t'\big)
\end{equation}
for $t\neq t'$ and varying $\mathbf{x}$. By the conditional independence of $U_n$ and $V_n$,
\begin{equation}\label{EqualityOfLaws}
\begin{array}{ccl}
{\cal L}\big(R_n\big|\mathbf{X}_n=\mathbf{x},T_n=t\big)\!\!&= &\!\!
{\cal L}\big(\rho(V_n,\mathbf{X}_n,T_n)\big|\mathbf{X}_n=\mathbf{x},T_n=t\big)\\
\!\!&= &\!\!
{\cal L}\big(\rho(V_n,\mathbf{x},t)\big|\mathbf{X}_n=\mathbf{x},\tau(U_n,\mathbf{x})=t\big)\\
\!\!&= &\!\!
{\cal L}\big(\rho(V_n,\mathbf{x},t)\big|\mathbf{X}_n=\mathbf{x}\big)\\
\!\!&= &\!\!
{\cal L}\big(\rho(V_n,\mathbf{x},t)\big),
\end{array}
\end{equation}
so any difference between the two laws of (\ref{ComparisonOfLawsI}) is purely a result of the difference between taking $t$ and $t'$ in the last argument of the response function,
and it follows that the function $t\!\rightarrow\! {\cal L}\big(\rho(V_n,\mathbf{X}_n,t)\big|\mathbf{X}_n\!\!=\!\mathbf{x}\big)$---which then provides a
full description of the causal effect of the treatment on the response under a situation $\mathbf{x}$---can be estimated from observed 
data in the guise of ${\cal L}\big(R_n\big|\mathbf{X}_n=\mathbf{x},T_n=t\big)$.

From this last law one can compute (and estimate from data) probabilities such as
\begin{equation}\label{FunctionalOfLawI}
\mathbf{P}\big(\rho(V_n,\mathbf{X}_n,t)\le r\big)=
\sum\nolimits_{\mathbf{x}}\mathbf{P}\big(\rho(V_n,\mathbf{x},t)\le r\big|\mathbf{X}_n=\mathbf{x}\big)\mathbf{P}(\mathbf{X}_n=\mathbf{x}),
\end{equation}
which determine the law of the potential outcomes, and quantities such as
\begin{equation}\label{FunctionalOfLawII}
\sum\nolimits_{\mathbf{x}}\left\{\mathbf{E}\big(\rho(V_n,\mathbf{x},t)\big|\mathbf{X}_n=\mathbf{x}\big)-
\mathbf{E}\big(\rho(V_n,\mathbf{x},t')\big|\mathbf{X}_n=\mathbf{x}\big)\right\}
\mathbf{P}(\mathbf{X}_n=\mathbf{x}),
\end{equation}
the overall mean difference between the responses to treatments $t$ and $t'$, which is
$\mathbf{E}\big[\rho(V_n,\mathbf{X}_n,t)-\rho(V_n,\mathbf{X}_n,t')\big]$---the expected difference between the 
potential outcomes to $t$ and $t'$ of \emph{the same} arbitrary individual.

\section{More detailed models; levels of causality}\label{LevelsOfCausality}

We begin by introducing some terminology\footnote{Proposed by Cox and Wermuth \cite{CoxWermuth04}.} which qualifies the assumptions
involved in a causal model and the degree to which those assumptions are thought to hold in a real-life problem. We speak of \emph{zero-level causality}
when a model like (\ref{CausalModelII}) is assumed \emph{without} unconfoundedness, so that conclusions based on samples of the $(\mathbf{X}_n,T_n,R_n)$s
may say something about the association between the treatment and the response but not necessarily about the causal effect of the former on the latter
(since at least some confounding subsists). It seems that observational studies are mostly carried out at the zero-level, even if they do not explicitly invoke model
(\ref{CausalModelII}) or do not discuss the extent of confounding, and even if they vary with respect to the statistical methods used.\footnote{In fact, many studies assume (at least implicitly) that the data follow a regression model, a very particular version of model (\ref{CausalModelII}) with unconfoundedness (see remark (i) following (\ref{CausalModelII})); sometimes, however, some of the actual data represented by
${\mathbf X}_n$ or $T_n$ are preceded by, or at least are synchronous with, those represented by $R_n$, which not only invalidates unconfounding as it perverts the notion of causality.}

In \emph{first-level causality} one assumes a model like (\ref{CausalModelII}) with
unconfoundedness---perhaps in a less explicit form, such as
Rubin's, discussed in remark (ii),
p.~\!\pageref{FootnoteUnconfoundednessAccordingToRubin}---and knows or believes that the model corresponds sufficiently well to the real-life
problem. In particular, one has identified most of the factors influencing both treatment and response---the confounders in $\mathbf{X}_n$---and in principle can
collect data on them in order to apply methods of stratification and matching to test for or estimate the causal effect of the treatment.
Indeed, we have seen that under unconfoundedness the effect of the treatment on the response in model (\ref{CausalModelII}), which is
completely characterized by ${\cal L}\big(\rho(V_n,\mathbf{x},t)\big)$ for varying $(t,\mathbf{x})\in{\cal T}\times {\cal X}$, can be estimated from observed data in the guise of ${\cal L}\big(R_n\big|\mathbf{X}_n\!=\!\mathbf{x},T_n=t\big)$.
However, as suggested above, not many studies are carried out at the first-level, and few present a thorough argumentation to justify unconfoundedness
based on a list of variables identified as confounders.  

Identifying which variables must be conditioned upon in order to make the assumption of unconfoundedness tenable is recognized by many (see for example section 3.7 of
\cite{ImbensRubin:2015}, pp.~\!76--7 of \cite{Rosenbaum2002}) as a difficult task, demanding substantial extra-statistical knowledge about the
real-life problem in hand.
As a general rule, researchers such as Rosenbaum and Rubin recommend that ``all the relevant covariates, that is, all the variables that may be associated with both
outcomes and assignment to treatment'' be considered as possible confounders, and say that ``there is little or no reason to avoid adjustment for a true covariate,
a variable describing subjects before treatment'' (\cite{ImbensRubin:2015},\cite{Rosenbaum2002}). Pearl and others (e.g.~\!\cite{Shrier2008},\cite{Pearl2009c}),
on the other hand, have pointed out that in certain situations conditioning on some variables may create rather than remove confounding (a phenomenon that will be illustrated below), so that in order to
justify conditioning on a particular set of variables one generally needs to postulate (and defend, based on extra-statistical knowledge) a \emph{detailed
causal model} describing the relationships between the various components of the situation, the treatment and the response---not just the model we have
considered in section \ref{StatisticalModel} of the cause and effect relationship between treatment and response given a situation.

We speak of \emph{second-level causality} whenever such a detailed causal model of the situation, treatment and response is assumed or can be
justified in a real-life problem. Since, as noted in section \ref{DeterministicModel}, a given problem may admit various levels of causality, second-level causality may be based on
several models of varying complexity (which of course should be consistent with each other), and, as we shall see below, it may be reduced to first-level causality.   

As an example of the `second-level model' we have in mind let us take
\begin{equation}\label{ExampleSecondLevelModel}
\left\{
\begin{array}{l}
X_1=\varphi_1(U_1), X_2=\varphi_2(U_2), \\
X_3=\varphi_3(U_3,X_1,X_2),\\
X_4=\varphi_4(U_4,X_1),\\
X_5=\varphi_5(U_5,X_2),\\
T=\tau(U,X_3,X_4),\\
X_6=\varphi_6(U_6,T),\\
R=\rho(V,X_3,X_5,X_6),
\end{array}
\right.
\end{equation} 
where $U$, $V$, $U_1,\ldots,U_6$ are independent standard uniform random variables, to be regarded as unobservable, purely random `errors' and therefore usually referred to as {\it exogenous},
and $\tau$, $\rho$, $\varphi_1,\ldots,\varphi_6$
are real-valued functions.\footnote{This model is used as an example in \cite{Pearl2009b},
pp.~\!113-4, and in \cite{PearlEtAl2016}, p.~\!64. In real-life problems, different 
variables in a set of equations like (\ref{ExampleSecondLevelModel}) may represent measurements of the same quantity made at different times, the
more recent one being a function of the earlier one.}

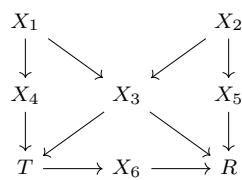
\begin{figure}[h!]
\begin{tikzcd}
& X_1 \arrow[d] \arrow[rd] &
& X_2 \arrow[d] \arrow[ld]\\
& X_4 \arrow[d]
& X_3 \arrow[rd] \arrow[ld]
& X_5 \arrow[d]\\
& T \arrow[r]
& X_6 \arrow[r]
& R
\end{tikzcd}
\caption{Graph of the model at the second-level of causality defined by (\ref{ExampleSecondLevelModel}).}\label{Figure1}
\end{figure}

The schematic representation of this set of equations in
figure \ref{Figure1}---a `causal graph'---omits the essential information about the exogenous random variables, but provides a more digestible summary of the
essential aspects of it, namely of which variables `influence' which. We shall
use standard or self-explanatory terminology when referring to such graphs: the variables $X_1,X_2,\ldots,T, R$ correspond to \emph{nodes} of the graph;
the graph is \emph{directed} because its edges are arrows; $X_1\rightarrow X_4\rightarrow T\rightarrow X_6 \rightarrow R$ is a \emph{directed path} leading $X_1$ to $R$, and we say that $X_1$ {\it leads to} $R$; $X_1\rightarrow X_4\rightarrow T\leftarrow X_3\rightarrow R$
is a \emph{path} (not a directed one); $T$ is a \emph{descendant} of $X_1$, and a \emph{direct descendant} of $X_4$ (and of $X_3$), and $X_1$ and $X_4$ are \emph{ancestors} of $T$ (and of $X_6$ and $R$); $T \leftarrow X_4$ is a subpath of the path $T \leftarrow X_4\leftarrow X_1$; and so on.
Note that a causal graph is {\it acyclic}: no two nodes point arrows to each other, and, more generally, a directed path from a node cannot lead back to the node (a future event cannot influence a past event). Actually, whenever we talk of a path we really mean a {\it simple path}: in it, a node appears only once (so even among undirected paths we never consider closed ones---those that begin and end in the same node).

\begin{remark}
We do not regard the graph of figure \ref{Figure1} as
a model: the model is the mapping that assigns a set or vector of numbers to each $\omega\in\Omega$, and the figure is a summary of this mapping which can serve as a tool to analyse it
in some respects. While in a general discussion one may consider figure \ref{Figure1} instead of (\ref{ExampleSecondLevelModel}), in a more specific example one may wish,
for example, that $\tau$ be constant in the second argument, in which case
the arrow from $X_3$ to $T$ may be misleading.\,\,$\square$
\end{remark}

Model (\ref{ExampleSecondLevelModel}) represents second-level causality because it consists of more than just the two equations of (\ref{CausalModelII}); but it is not the only such model
since the model without the equation $X_6=\varphi_6(U_6,T)$ and with $X_6$ replaced by
$\varphi_6(U_6,T)$ in the definition of $R$---a graph of which
is obtained by deleting the node $X_6$ from figure \ref{Figure1} and linking $T$ to $R$ by a single arrow---is simpler but still of second-level causality. Evidently, if we discard $X_6$ and put
$\mathbf{X}=(X_1,\ldots,X_5)$ then the model satisfies (\ref{CausalModelII}) with unconfoundedness and in essence reduces to first-level-causality.
More generally, if some of the functions and variables representing `intermediate relations' are
substituted into the functions they actuate, one may reduce the model to a
first-level model: For instance, from (\ref{ExampleSecondLevelModel}) follow
\begin{equation*}
R=\rho(V,X_3,X_5,X_6)=\rho\big(V,X_3,\varphi_5(U_5,X_2),\varphi_6(U_6,T)\big)
\,=:\tilde\rho\big(\tilde V,\mathbf{X},T\big),
\end{equation*}
\[
T=\tau(U,X_3,X_4)=\tau\big(U,X_3,\varphi_4(U_4,X_1)\big)=:\tilde\tau\big(\tilde U,\mathbf{X}\big),\quad
\]
where $\mathbf{X}=(X_1,X_2,X_3)$ and $\tilde U$ and $\tilde V$ are standard uniforms obtained from $(U_4,U)$ and $(U_5,U_6,V)$, respectively,
independent and independent of $\mathbf{X}$, so the conditions of the basic model (\ref{CausalModelII}) are satisfied with unconfoundedness
and\footnote{$\tilde U$ and $\tilde V$ are independent and independent of $\mathbf{X}$ because the event $\{\mathbf{X}=\mathbf{x}\}$ is equal to 
$\big\{\big(\varphi_1(U_1),\varphi_2(U_2),\varphi_3\big(U_3,\varphi_1(U_1),\varphi_2(U_2)\big)\big)=\mathbf{x}\big\}$ and hence involves $(U_1,U_2,U_3)$,
while $\tilde U$ involves $(U_4,U)$ and $\tilde V$ involves $(U_5,U_6,V)$. In particular, the conditional distribution of $\tilde V$ given $\{\mathbf{X}=\mathbf{x}\}$
is independent of $\mathbf{x}$, whence the last equality in (\ref{ExampleSecondLevelModelB})
(cf. the remark on p.~\pageref{FootnoteConditioningCanBeRemoved}).
As we shall see below, other reductions to a first-level model are possible.
Note that $\tilde V$ may, for example, be defined by arranging the digits
in the decimal expansions of $V$, $U_5$, $U_6$ alternately, so that $\tilde\rho$ involves first the `unfolding' of a sequence in three sequences of digits.}
\begin{equation}\label{ExampleSecondLevelModelB}
{\cal L}\big(R\big|\mathbf{X}=\mathbf{x},T=t\big)\!=\!
{\cal L}\big(\tilde\rho\big(\tilde V,\mathbf{x},t\big)\big|\mathbf{X}=\mathbf{x},\tilde\tau\big(\tilde U,\mathbf{x}\big)=t\big)
\!=\!{\cal L}\big(\tilde\rho\big(\tilde V,\mathbf{x},t\big)\big).
\end{equation}

Despite the difficulty of setting up a convincing second-level model in a real-life problem, it seems evident
that any earnest effort in that direction can only help in thinking about and defending the assumption of unconfoundedness in a corresponding first-level model---and  
in opening it to scrutiny and criticism.\footnote{Pearl has been explaining this for years (see the references in footnote \ref{ThreeReferencesToPearl}); for a different opinion see \cite{ImbensRubin1995} and p.~\!22 of \cite{ImbensRubin:2015}.
It must be admitted, however, that if in relatively simple problems, such as (for instance) those considered in \cite{RobinsWasserman1997}, setting up and justifying a second-level model is possible, in many observational studies of the type considered in \cite{Rosenbaum2002} and
\cite{ImbensRubin:2015} it may be simply impossible. In fields such as Economics and Sociology,
identifying and estimating causal effects seems to be particularly difficult; see the valuable analyses of Freedman in the later chapters of \cite{Freedman:2009} and in
chapter 5 of \cite{Freedman:2010}.}
The study of such a second-level model can also render a statistical analysis more efficient if it is found that unconfoundedness holds by conditioning on a smaller or more easily available subset of variables, although the actual measure of causal effect depends somewhat on the subset chosen: thus, essentially the same method of reduction just used with $\mathbf{X}=(X_1,X_2,X_3)$ shows that
both $\mathbf{X}=(X_1,X_3)$ and $\mathbf{X}=(X_2,X_3)$ guarantee unconfoundedness, 
but (\ref{ExampleSecondLevelModelB}) could be easier to estimate with the latter and a somewhat different measure of causal effect could be estimated
with the former.\footnote{\label{FootnoteVerificationX2X3}Take for instance $\mathbf{X}=(X_2,X_3)$:
${\cal L}(R|T=t,\mathbf{X}=\mathbf{x}) 
={\cal L}\big(\rho(V,x_3,\varphi_5(U_5,x_2),\varphi_6(U_6,t))\big|$ 
$\tau(U,x_3,\varphi_4(U_4,X_1))=t,\mathbf{X}=(x_2,x_3)\big)=
{\cal L}\big(\rho(V,x_3,\varphi_5(U_5,x_2),\varphi_6(U_6,t))\big)$. It is important to remember that different reductions to first-level causality usually lead to different measures of causal effect; this will be made clear in the remark of subsection \ref{PearlCriterion}.} 

Besides, it seems that the possibility mentioned above that conditioning on certain variables may create rather than remove confounding, although considered
by some to be mostly theoretical (e.g.~\!\cite{Rubin2009}), cannot be ruled out. For instance, if we try to reduce model (\ref{ExampleSecondLevelModel}) to
first-level causality by conditioning on $\mathbf{X}:=X_3$ alone, writing
\[
R=\rho(V,X_3,X_5,X_6)=\rho\big(V,X_3,\varphi_5(U_5,X_2),\varphi_6(U_6,T)\big)
=:\check\rho\big(\check V,\mathbf{X},T\big),
\]
\[
T=\tau(U,X_3,X_4)=\tau\big(U,X_3,\varphi_4(U_4,X_1)\big)=:\check\tau\big(\check U,\mathbf{X}\big),\quad
\]
we find that $\check U$ and $\check V$ are no longer conditionally independent because
 $\check U$ involves $(U_4,U,X_1)$, $\check V$ involves $(U_5,U_6,V,X_2)$, and
$\check U$ and $\check V$ are `entangled' by $X_1$ and $X_2$ through the conditioning event
$\{\mathbf{X}\!=\!\mathbf{x}\}\!=\!\{\varphi_3\big(U_3,X_1,X_2\big)\!=\!\mathbf{x}\}$, and conclude that
unconfoundedness fails if $X_3$ is the only variable to be corrected for in (\ref{ExampleSecondLevelModel}); but if $\rho$ and $\tau$ are both constant in their second arguments (so that neither $R$ nor $T$ involve $X_3$) then (\ref{ExampleSecondLevelModel}) is automatically a first-level model whose unconfoundedness is destroyed by conditioning on $X_3$.\footnote{For another example see pp.~\!118-20 of \cite{Ferreira:2015}.}

Although a model at the second-level of causality ought to be based principally on extra-statistical knowledge, it has been pointed out by Pearl and others
(e.g.~section 2.5 of \cite{PearlEtAl2016}) that such a model implies the conditional independence of certain variables and hence can be checked, in part, by means of data on those variables.
Thus, in model (\ref{ExampleSecondLevelModel}) $X_3$ and $X_4$ are independent conditionally on $X_1$ because
\begin{eqnarray*}
{\cal L}\big(X_4\big|X_1=x_1,X_3=x_3\big)\!\!\!&=&\!\!\!{\cal L}\big(\varphi_4(U_4,x_1)\big|X_1=x_1,\varphi_3(U_3,x_1,X_2)=x_3\big)\\
\!\!\!&=&\!\!\! {\cal L}\big(\varphi_4(U_4,x_1)\big)
\end{eqnarray*}
does not depend on $x_3$, so if data are available on $X_1$, $X_3$ and $X_4$ a test of independence provides a test of the model;
similarly, the independence of $X_3$ and $X_5$ conditionally on $X_2$, the independence of $X_1$ and $X_6$ conditionally on $T$, etc.,
are implications of the model that can be used to criticize it. On the other hand,
as in other questions of
goodness-of-fit, unless the sample size is very large the non-rejection of conditional independence does not imply the approximate correctness of the portion of the model being 
tested.

\subsection{Reduction to first-level causality; Pearl's criterion}\label{PearlCriterion}

Let us think of a generic model at the second-level of causality involving a number of random variables $X_1,X_2,\ldots,T,R$, such that $T$ and $R$ are functions of some of the $X_j$s,
none of the $X_j$s is a function of $(T,R)$ and $R$ is a function of $T$, and consider the problem of finding a set $\mathbf{X}$ of $X_j$s in the model which allows the estimation of the effect of $T$ on $R$ by conditioning upon it---a set  which provides the means of  reducing the model to a first-level model.
There exists at least one set with the required property, namely the union
$\mathbf{X}\!=\!\mathbf{X}_1\cup\mathbf{X}_2$ of the set $\mathbf{X}_1$ of $X_j$s that affect $T$ with the set $\mathbf{X}_2$ of $X_j$s that affect $R$ (and which by our definitions are not affected by $T$);
for then $R=\rho(V,\mathbf{X}_2,T)$, $T=\tau(U,\mathbf{X}_1)$, $\mathbf{X}$ does not involve $T$, and so we have, with $\mathbf{x}$, $\mathbf{x}_1$ and $\mathbf{x}_2$ numerical vectors such that
$\mathbf{x}=(\mathbf{x}_1,\mathbf{x}_2)$ and $\mathbf{X}=\mathbf{x}$
$\Leftrightarrow$ $\mathbf{X}_1=\mathbf{x}_1$ $\wedge$ $\mathbf{X}_2=\mathbf{x}_2$, 
\begin{equation}\label{EquationSelection0}
{\cal L}(R|T=t,\mathbf{X}=\mathbf{x})=
{\cal L}\big(\rho(V,\mathbf{x}_2,t)\big|
\tau(U,\mathbf{x}_1)=t, \mathbf{X}=\mathbf{x}\big)=
{\cal L}\big(\rho(V,\mathbf{x}_2,t)\big)
\end{equation}
by the independence of $U$, $V\!$ and the exogenous variables involved in $\mathbf{X}$, and the last term here characterizes the causal effect of $T$ on $R$ completely and can be estimated from observed data in the guise of the first term. The point is that there may be a smaller, or in some sense more convenient, $\mathbf{X}$.\footnote{In the complementary notes at the end of this subsection we consider situations in which some of the $X_j$s are functions of $T$. 
Note that if some of the $X_j$s are functions of $T$ and are also involved in $R$ then we can write them explicitly in terms of $T$ and other variables in order to achieve a `reduced' expression for $R$ that satisfies our requirements. For instance, in model (\ref{ExampleSecondLevelModel}) the response can be written as
$R=\rho\big(V,X_3,X_5,\varphi_6(U_6,T)\big)\equiv\tilde\rho\big(\tilde V\!,X_3,X_5,T\big)$ 
so to exclude $X_6$ from $\mathbf{X}_2$.}

\begin{figure}[h!]
\begin{tikzcd}
& \mathbf{\check X} \arrow[r, dashrightarrow] \arrow[r, dashleftarrow]
\arrow[dd]  \arrow[rr, dashrightarrow, bend left]\arrow[rr, dashleftarrow, bend left]
\arrow[rrdd, dashrightarrow] 
& \mathbf{X} 
\arrow[rdd,dashrightarrow] \arrow[ldd,dashrightarrow]
\arrow[r, dashrightarrow] \arrow[r, dashleftarrow]
&   \mathbf{\hat X}\arrow[dd] \arrow[lldd,dashrightarrow] \\
&\,&\,&\,\\
& T \arrow[rr]  & \,  & R
\end{tikzcd}
\caption{Graphical representation of a second-level model to which Pearl's criterion applies.}\label{IllustrationChoiceOfVariables}
\end{figure}
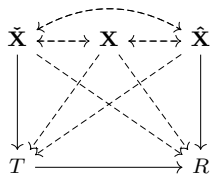

In what follows, denote by $\mathbf{X}$ a (possibly empty) set with the required property, or a candidate for such a set. 
Let $\mathbf{\check X}$ be the set of variables not in $\mathbf{X}$ that affect the treatment as arguments of $\tau$ and $\mathbf{\hat X}$ the set of variables not in $\mathbf{X}$ that affect the response as arguments of $\rho$. 
Figure \ref{IllustrationChoiceOfVariables} sketches a graph of the model
with the sets $\mathbf{X}$, $\mathbf{\check X}$ and $\mathbf{\hat X}$, and with the double, dashed arrows indicating that 
some nodes $X_i\!\in\!\mathbf{X}$, $\check X_j\!\in\!\mathbf{\check X}$ and $\hat X_k\!\in\!\mathbf{\hat X}$ may be functions of each other---so that we may have, for instance,
$X_i=\varphi_i(U_i,\check X_j)$, $\check X_j=\check\varphi_j({\check U}_j,\hat X_k)$,
which in a detailed representation of the graph would correspond to the path
$X_i \leftarrow \check X_j \leftarrow \hat X_k$.\footnote{As always, the $U_i$s, $\check U_j$s and $\hat U_k$s are independent standard uniform random variables. Note that if the equation
$X_i=\varphi_i(U_i,\check X_j)$ is used to define a second-level model then that same model cannot
be defined with an equation $\check X_j=\varphi_j(U_i,X_i)$ for some $\varphi_j$ (in terms of the graph this means that each node may point an arrow to a node or have an arrow pointed to it but no two nodes can point arrows to each other), although the first equation can in principle be inverted, globally or locally, to give $\check X_j$ as a function $\varphi_i^{-1}$ of
$X_i$ and $U_i$; consider for instance the model $X_1=U_1$,
$X_2=U_2 X_1 $ (cf.~\!footnote \ref{FootnoteConcern} and the text leading to it).} 

Since $\tau$ and $\rho$ generally have common arguments,
$\mathbf{\check X}$ and $\mathbf{\hat X}$ are typically not disjoint, but by definition none of their elements is to be conditioned upon and hence $\mathbf{X}\cap\mathbf{\hat X}=\mathbf{X}\cap\mathbf{\check X} =\emptyset$. Since not all variables need to be involved in $\tau$ and $\rho$, $\mathbf{X}\cup\mathbf{\check X}\cup\mathbf{\hat X}$ need not contain all the $X_j$s in the model; thus figure \ref{IllustrationChoiceOfVariables} may correspond to a somewhat incomplete representation of the model.

To help fixing ideas, consider conditioning on $\mathbf{X}\!=\!\{X_2,X_3\}$ in model (\ref{ExampleSecondLevelModel}) without the equation defining $X_6$ and with $R\!=\!\rho\big(V,X_3,X_5,\varphi_6(U_6,T)\big)$.
Since $T\!=\!\tau(U,X_3,X_4)$, we may take $\mathbf{\check X}\!=\!\{X_4\}$
and $\mathbf{\hat X}\!=\!\{X_5\}$.
But because we can also write
$T\!=\!\tau\big(U,X_3,\varphi_4(U_4,X_1)\big)$, we may take 
$\mathbf{\check X}\!=\!\{X_1\}$ and $\mathbf{\hat X}\!=\!\{X_5\}$ instead;
and because we can also write $R\!=\!\rho\big(V,X_3,\varphi_5(U_5,X_2),\varphi(U_6,T)\big)$
we may take $\mathbf{\check X}\!=\!\{X_1\}$ and $\mathbf{\hat X}\!=\!\emptyset$ instead.

If $\mathbf{X}$ is to allow the estimation of the effect of $T$ on $R$
we must have\footnote{Although $\mathbf{X}\cap\mathbf{\hat X}=\emptyset$, in general we are not allowed to drop the conditioning on $\{\mathbf{X}=\mathbf{x}\}$ from the last term of
(\ref{EquationSelection1}) because $\mathbf{X}$ and $\mathbf{\hat X}$ may be dependent.}
\[
{\cal L}(R|T=t,\mathbf{X}=\mathbf{x})=
{\cal L}\big(\rho(V,\mathbf{\hat X},\mathbf{x},t)\big|
\tau(U,\mathbf{\check X},\mathbf{x})=t, \mathbf{X}=\mathbf{x}\big)\vspace{-0.5cm}
\]
\begin{equation}\label{EquationSelection1}
\quad\vspace{-0.25cm}
\end{equation}
\[
\quad\,\,\,\,\,\,\,={\cal L}\big(\rho(V,\mathbf{\hat X},\mathbf{x},t)
\big|\mathbf{X}=\mathbf{x}\big).
\]
But the second equality holds if and only if $\mathbf{\check X}$ and $\mathbf{\hat X}$ are
independent conditionally on $\mathbf{X}$. Indeed, if this condition holds then
 $\tau(U,\mathbf{\check X},\mathbf{x})\!=\!t$ can be dropped from the second term in (\ref{EquationSelection1});
but if it fails then the treatment may, due to the constraint $\tau(U,\mathbf{\check X},\mathbf{x})\!=\!t$, exert an effect on the response through $\mathbf{\hat X}$---not only through the last argument of $\rho$.
The condition stated characterizes the desired sets $\mathbf{X}$, but it does 
not yet provide a direct, workable means of identifying them.

For the second equality in (\ref{EquationSelection1}) to fail there must be at least two variables $\check X\in\mathbf{\check X}$ and $\hat X\in\mathbf{\hat X}$ which are dependent given $\mathbf{X}$.
In the graph of a second-level model there is always a path (not necessarily a directed one)
between nodes, certainly if they are dependent. Thus the dependence between $\check X$ and $\hat X$
implies that there is a path between them which links $T$ and $R$
(because $\check X$ points an arrow to $T$ and $\hat X$ points one to $R$) and which represents a chain of random equations involving $\check X$ and $\hat X$ (hence also $T$ and $R$) and other nodes.

As an example, and in order to introduce some notation used below, let us consider the path
\begin{equation}\label{EquationSelection2}
T\leftarrow\check X \leftarrow N_1 \rightarrow N_2 \leftarrow N_3 \rightarrow N_4
\leftarrow N_5 \rightarrow \hat X \rightarrow R
\end{equation}
with $N_1,N_2,\ldots$ standing for generic nodes (not necessarily belonging to $\mathbf{X}$, $\mathbf{\check X}$ or $\mathbf{\hat X}$). This path represents a chain
of equations of the type
\[
T=\tau(U,\check X,\ldots),\quad
\check X=\check\varphi(\check U,N_1,\ldots),\quad
N_2=\psi_2(W_2,N_1,N_3,\ldots),
\]
\[
N_4=\psi_4(W_4,N_3,N_5,\ldots),\quad
\hat X=\hat\varphi(\hat U,N_5,\ldots),\quad
R=\rho(V,\hat X,\ldots,T),
\]
where $U,\check U,W_2,...,V$ are independent standard uniforms,
$\check\varphi, \hat\varphi, \psi_2,\ldots$ are certain functions, and the ellipsis after a list of arguments indicates the possible presence of other variables not in (\ref{EquationSelection2}).

Pearl proposed what appears to be the most general criterion for deciding which variables to include in
$\mathbf{X}$ (or not to include, when $\mathbf{X}$ is empty): If $\mathbf{X}$ has the following property, then it secures the conditional independence of $\mathbf{\check X}$ and $\mathbf{\hat X}$ and hence the second identity in (\ref{EquationSelection1}):\footnote{It is due to Pearl;
see pp.~\!17--8, 79--81 of \cite{Pearl2009a},
pp.~\!106, 114 of \cite{Pearl2009b}, and the earlier references provided in these sources.
We have not been satisfied with any proof of this result, but have not yet studied 
the more recent proof in \cite{Perkovic2018}.
A computer program is almost a necessity if one is to apply the criterion to large models;
for flexible software see \cite{TextorEtAl2016}.} 

\vspace{0.25cm}
$\mathbf{X}$ is such that for each path $P$ linking $T$ to $R$ with an arrow pointing at $T$
and an arrow pointing at $R$ one (hence only one) of these conditions holds:
\begin{description}
\item [{\bf (i)}\,] There is a node $X\in\mathbf{X}$ in $P$ pointing an arrow to a node in $P$, so $P$ contains a subpath of type $N'\rightarrow X\rightarrow N''$ or of type
$N'\leftarrow X\rightarrow N''$;
\item [{\bf (ii)}] No node of $\mathbf{X}$ in $P$ points an arrow to a node in $P$,
and $P$ contains a subpath of the form $N'\rightarrow N\leftarrow N''$ where neither $N$ nor any of its descendants is in $\mathbf{X}$.
\end{description}
\vspace{0.1cm}

Let us illustrate the application of the criterion to the model of figure \ref{Figure1} without the node $X_6$:
The path $T\leftarrow X_3\rightarrow R$ does not permit
{\bf (ii)}, so {\bf (i)} must hold and hence $X_3$ must be included in $\mathbf{X}$.
But if {\bf (i)} is to hold with the path
\[
T\leftarrow X_4\leftarrow X_1 \rightarrow X_3 \leftarrow X_2
\rightarrow X_5 \rightarrow R
\]
another node must be added to $\mathbf{X}$ (for otherwise $X_3$,
the single node in it, would point no arrows in that path), and this may be $X_1$, $X_2$,
$X_4$ or $X_5$. Since $X_3$ caters for the remaining two paths linking $T$ to $R$, we conclude that
 if $\mathbf{X}$ consists of $X_3$ and at least one of $X_1$, $X_2$,

\section{Specification and identification of causal effects: some examples}\label{IdentificationEstimation}

Our purpose in this section is to study a number of `causal problems' that have been considered (as examples, mostly) in the literature, and to solve them by elementary probability arguments.

\subsection{Smoking and the genotype theory}\label{SectionSmokingAndGenotype}

Pearl (pp.~\!83--84 of \cite{Pearl2009a}) considers the model
\begin{equation}\label{ExampleSmokingAndGenotype}
\left\{
\begin{array}{l}
X=\varphi_0(U_0), \\
Y=\varphi_1(X,U_1), \\
Z=\varphi_2(Y,U_2),\\
W=\varphi_3(X,Z,U_3),
\end{array}
\right.
\end{equation}
where $U_0$, $U_1$, $U_2$, $U_3$ are independent standard uniforms,
$W$, $Y$, $Z$ are regarded as observable,
$X$ and the $U_i$s as unobservable, and the $\varphi_i$s as unknown real-valued functions.
The model, summarized by figure \ref{FigureSmokingAndGenotype}, is intended to describe the combined effects on the development of lung cancer of smoking and of a putative genotype that is not only carcinogenic but also predisposes people to smoke;
the random variables pertain to an individual randomly drawn from some population, $X$ stands for the individual's genotype, $Y$ for his level of smoking and $Z$ for the concentration of tar deposits in his lungs determined at some point in time, and $W$ for the indicator of whether the individual develops lung cancer at a later time.

\begin{figure}[h!]
\begin{tikzcd}
& 
X \arrow[rd] \arrow[ld]
\\
 Y \arrow[r] & Z \arrow[r] & W
\end{tikzcd}
\caption{Graph of model (\ref{ExampleSmokingAndGenotype}).}\label{FigureSmokingAndGenotype}
\end{figure}

The question of interest is whether the effect of smoking on the development of lung cancer
(the effect of $Y$ on $W$) can be estimated despite the confounder of that effect ($X$) being unobservable; the answer is that it can, thanks to the knowledge about the concentration of tar in the lungs ($Z$), a more immediate cause of cancer.
This model is probably not very realistic, but the question behind it is not trivial and could be of interest in similar problems.\footnote{See Freedman's critique of the model on
pp.~\!272-3 of \cite{Freedman:2010}. In Freedman's opinion, the `perverse theory' about smoking and lung cancer, according to which it is the genotype rather than smoking that causes cancer, is refuted by carefully designed population studies, not by statistical analyses of `detailed' models such as (\ref{ExampleSmokingAndGenotype}).
As far as we know, Freedman (pp.~\!267--9 of \cite{Freedman:2010}) was the first to prove Pearl's result by means of
standard probability arguments; our approach is less general 
because it assumes the variables to be discrete, but we make an explicit connection
between the quantity to be estimated from data on $(W,Y,Z)$ and the treatment
effect---comparing the proofs is easy since we follow Freedman's notation rather than Pearl's.}

Model (\ref{ExampleSmokingAndGenotype}) implies 
\begin{equation}\label{ExampleSmokingAndGenotype0}
W=\varphi_3\left(X,\varphi_2(Y,U_2),U_3\right),\quad
Y=\varphi_1(X,U_1),
\end{equation}
which is of type (\ref{CausalModelII}) with unconfoundedness: the `treatment' $Y$ is a function of $X$, the response $W$ is a function of the treatment $Y$ and of $X$, the latter confounds the effect of $Y$ on $W$, but since $U_1$ and $U_2$ are independent (conditionally on $X$
as well as unconditionally) the treatment effect
is fully characterized by ${\cal L}(W|X=x,Y=y)$ and may be quantified through functionals of it, such as
\begin{equation}\label{ExampleSmokingAndGenotype1}
p_y(w):=\sum\nolimits_x \mathbf{P}(W=w|X=x,Y=y)\mathbf{P}(X=x).
\end{equation}
This can be estimated from data on $(W,X,Y)$, but our assumption is that only data on
$(W,Y,Z)$ are available. In order to show that one can write $p_y(w)$ in terms of probabilities pertaining to $W$, $Y$ and $Z$, we first put
\[
m_z(w):=\sum_x \mathbf{P}(W=w|X=x,Z=z)\mathbf{P}(X=x)\quad\quad\quad\quad\quad
\quad\quad\quad\quad\quad\quad\quad\,\,
\]
\vspace{-0.17cm}
\[
\quad\quad\,\,=\sum_x \mathbf{P}\!\left[\varphi_3(x,z,U_3)\!=\!w|\varphi_0(U_0)\!=\!x,\varphi_2(\varphi_1(x,U_1),U_2)\!=\!z\right]\mathbf{P}(X\!=\!x)
\]
\[
=\sum_x \mathbf{P}(\varphi_3(x,z,U_3)=w)\mathbf{P}(X=x),
\quad\quad\quad\quad\quad\quad\quad\quad\quad\quad\quad
\]
where in the last equality we make use of the independence of the $U_i$s, and
\begin{equation}\label{ExampleSmokingAndGenotype00}
l_y(w):=\sum\nolimits_z m_z(w)\mathbf{P}(Z=z|Y=y),
\end{equation}
and then note that
\[
p_y(w)=\sum_x \mathbf{P}
\left[\varphi_3(x,\varphi_2(y,U_2),U_3)\!=\!w\left|\varphi_0(U_0)\!=\!x,\varphi_1(x,U_1)\!=\!y\right.\right]
\mathbf{P}(X\!=\!x)
\]
\vspace{-0.17cm}
\[
\quad\,\,=\sum_x\mathbf{P}\left[\varphi_3(x,\varphi_2(y,U_2),U_3)=w\right]\mathbf{P}(X=x)
\quad\quad\quad\quad\quad\quad\quad\quad\quad
\]
\[
\quad\,\,=\sum_{x,z}\mathbf{P}\left[\varphi_3(x,\varphi_2(y,U_2),U_3)=w,
\varphi_2(y,U_2)=z\right]\mathbf{P}(X=x)\quad\quad\quad\!\!
\]
\[
\quad\,\,=\sum_{z}\mathbf{P}\left[\varphi_2(y,U_2)=z\right]
\sum_{x}\mathbf{P}\left[\varphi_3(x,z,U_3)=w\right]\mathbf{P}(X=x)\quad\quad\quad\,
\]
\[
\quad\,\,=\sum_{z}\mathbf{P}\left[\varphi_2(y,U_2)=z\right]
m_z(w)
=\sum_{z}\mathbf{P}(Z=z|Y=y)\,m_z(w)\quad\quad
\]
\[
\quad\,\,=l_y(w),\quad\quad\quad\quad\quad\quad\quad\quad\quad\quad\quad
\quad\quad\quad\quad\quad\quad\quad\quad\quad\quad\quad\quad\quad\,\,\,
\]
the penultimate step following from
\[
\mathbf{P}(Z=z|Y=y)=\mathbf{P}[\varphi_2(y,U_2)=z|\varphi_1(\varphi_0(U_0),U_1)=y]=
\mathbf{P}(\varphi_2(y,U_2)=z).
\]
Thus the result will be proved if we prove that $m_z(w)$ can be written in terms
of probabilities pertaining to $W$, $Y$ and $Z$.

By the independence of $X$, $U_1$, $U_2$ and $U_3$, we have
\[
\mathbf{P}(W=w|Y=y,Z=z)=
\frac{
\mathbf{P}\left[\varphi_3(X,z,U_3)=w,\varphi_1(X,U_1)=y,\varphi_2(y,U_2)=z\right]
}
{
\mathbf{P}(Y=y,Z=z)
}
\]
\[
\quad\quad\quad\quad\,\,\,\,=
\sum_x\frac{
\mathbf{P}\!\left({\varphi_3(x,z,U_3)=w,\,\varphi_1(x,U_1)=y,}\atop {\varphi_2(y,U_2)=z,
X=x}\right)
}
{
\mathbf{P}(\varphi_1(X,U_1)=y,\varphi_2(y,U_2)=z)
}
\]
\[
\quad\quad\quad\quad\quad\quad\quad\quad\quad\quad\,\,\,\,=
\sum_x
\frac{
\mathbf{P}(\varphi_3(x,z,U_3)=w)
\mathbf{P}\!\left({\varphi_1(x,U_1)=y,\,\varphi_2(y,U_2)=z,}\atop {
X=x}\right)
}
{
\mathbf{P}(\varphi_1(X,U_1)=y)\mathbf{P}(\varphi_2(y,U_2)=z)
}
\]
\[
\quad\quad\quad\quad\quad\quad\,\,\,=
\sum_x
\frac{
\mathbf{P}(\varphi_3(x,z,U_3)=w)
\mathbf{P}\!\left({\varphi_1(X,U_1)=y,}\atop {
X=x}\right)
}
{
\mathbf{P}(\varphi_1(X,U_1)=y)
}
\]
\[
\quad\quad\quad\quad\quad\quad\,\,\,\,\,
=\sum_x\mathbf{P}(\varphi_3(x,z,U_3)=w)\mathbf{P}(X=x|Y=y)
\]
(incidentally, the last couple of steps show that $X$ and $Z$ are independent given $Y$).
Since
\[
\mathbf{P}(\varphi_3(x,z,U_3)=w)=
\mathbf{P}(W=w|Z=z,X=x)
\]
(an observation already used in connection with the definition of $m_z(w)$), the preceding identity is equivalent to
\[
\mathbf{P}(W=w|Y=y,Z=z)=\sum_x\mathbf{P}(W=w|Z=z,X=x)\mathbf{P}(X=x|Y=y).
\]
Finally, integrating both sides here with respect to the distribution of $Y$ we get
\[
\sum_{y'} \mathbf{P}(W\!=\!w|Y\!=\!y',Z\!=\!z)\mathbf{P}(Y\!=\!y')=
\sum_{x} \mathbf{P}(W\!=\!w|Z\!=\!z,X\!=\!x)\mathbf{P}(X\!=\!x)
\]
\[
\quad\quad\quad\quad\!\!
=m_z(w),
\]
and the first term involves only probabilities pertaining to $W$, $Y$ and $Z$.

In conclusion, the effect of $Y$ on $W$ in the form of (\ref{ExampleSmokingAndGenotype1}) can be estimated from data on $(W,Y,Z)$ because
\[
p_y(w)=\sum_{y'\!,z} \mathbf{P}(W\!=\!w|Y\!=\!y',Z\!=\!z)\mathbf{P}(Y\!=\!y')
\mathbf{P}(Z\!=\!z|Y\!=\!y)=l_y(w).
\]

This result has been obtained by identifying the effect of smoking on cancer through
the reduced model (\ref{ExampleSmokingAndGenotype0}) and using the information
of the full model (\ref{ExampleSmokingAndGenotype}) to derive an alternative
expression---that of $l_y (w)$ in (\ref{ExampleSmokingAndGenotype00})---for a measure
of that effect---namely $p_y (w)$.
It is worth noting, however, that $l_y (w)$ has its own causal interpretation
within the full causal mechanism:
Under (\ref{ExampleSmokingAndGenotype}), $m_z (1)$ 
quantifies the effect of a concentration $z$ of deposited tar on the development of lung cancer irrespectively of genotype, since, for example, $m_z(1)-m_{z'}(1)$ is equal to
\[
\sum\nolimits_x\left\{{\mathbf P}(W=1|X=x,Z=z)-
{\mathbf P}(W=1|X=x,Z=z')\right\}{\mathbf P}(X=x),
\]
the expected difference between the `risk' of cancer at different concentrations $z$ and $z'$ of deposited tar for an individual/genotype randomly drawn from the population.
The differences $m_z(1)-m_{z'}(1)$ are meaningful because under model (\ref{ExampleSmokingAndGenotype}) the development of lung cancer is purely a function of genotype and tar accumulation and of other, `exogenous factors' embodied by $U_3$. On the other hand, tar deposits occur only through smoking (in particular they are not affected by genotype), so ${\mathbf P}(Z=z|Y=y)$ is a proper measure of the effect of smoking on tar accumulation. Thus (see (\ref{ExampleSmokingAndGenotype00})), $l_y (1)$ quantifies the effect of smoking on cancer by averaging the risk of cancer due to tar accumulation with respect to the distribution of the tar accumulation that results from smoking at level $y$.

\subsection{Case-control studies}
\label{SectionCaseReferentStudy}

Let us consider model (\ref{CausalModelII}) with unconfoundedness, with the $R_n$s and $T_n$s binary and the ${\mathbf X_n}$s as well as the $(U_n,V_n)$s independent and identically distributed, so that the vectors $({\mathbf X_n},T_n,R_n)$ too are independent and identically distributed, and think of
$R_n$ as an individual's status of a certain disease (equal to 1 if the individual is diseased), of $T_n$ as the individual's indicator of exposure to an agent suspected of causing the disease (equal to 1 if the individual has been exposed), and of ${\mathbf X_n}$ as a set of background characteristics  (year of birth, sex, occupation, etc.) which may influence the individual's exposure and disease statuses.
For simplicity we assume that the ${\mathbf X_n}$s are discrete, so that the
$({\mathbf X_n},T_n,R_n)$s too are discrete. As usual, we write $({\mathbf X},T,R)$ for a random vector with the same distribution as the $({\mathbf X_n},T_n,R_n)$s.

In situations where the disease is rare it may be necessary to collect a very large random sample of individuals in order to obtain enough {\it cases}---i.e.~\!diseased individuals---to provide evidence of the supposed effect of the exposure on the disease.
Instead of collecting a random sample one may then think of obtaining a substantial number of cases from readily accessible subpopulations---typically patients in hospitals that treat the disease---and comparing them in terms of frequency of exposure with a more general subpopulation which, due to the rarity of the disease, will consist mostly of nondiseased individuals. The intuition behind this proposal is that if the background characteristics of the cases are very similar to those of the nondiseased individuals and if the exposure indeed contributes to the disease then a greater rate of exposure must be found among the cases than among the nondiseased. But for it to work it is evident that great care is needed to ensure that the subpopulation of cases and the more general subpopulation are comparable enough with respect to everything except for frequency of disease and possibly frequency of exposure.

In a case-control study one tries to make the subpopulations comparable by matching each case with a `control', namely an individual from a more general population who has the same or practically the same background characteristics as the case. In the language of our model, this type of study can be described as follows: First, conditionally on the set $\{R_1,R_2,\ldots\}$ of responses of the whole population of interest, a sequence of cases is obtained by selecting random indices $I_1<I_2<\cdots$ from $\left\{n: R_n=1\right\}$, and a sample of $N$ cases is obtained from the indices in ${\cal C}_N=\{I_1,I_2,\ldots,I_N\}$ by setting
\[
({\mathbf X_{2n-1}'},T'_{2n-1},R'_{2n-1})\equiv
({\mathbf X_{2n-1}'},T'_{2n-1},1):=
({\mathbf X_{I_n}},T_{I_n},1)
\]
for $n=1,2,\ldots,N$.  
Then, conditionally on $\{({\mathbf X_{1}},R_1),({\mathbf X_{2}},R_2),\ldots\}$  (the pairs of background characteristics and responses of the whole population), for each ${\mathbf x}\in {\cal X}$ 
a sequence of random indices $J_{1}({\mathbf x})<J_{2}({\mathbf x})<\cdots$ is selected from $\left\{n: {\mathbf X}_n={\mathbf x}\right\}\setminus {\cal C}_N$
independently of all the other variables
and a corresponding sample of controls is obtained as
\[
({\mathbf X_{2n}'},T'_{2n},R'_{2n}):=
({\mathbf X_{2n-1}'},T_{J_n({\mathbf X_{2n-1}')}},R_{J_n({\mathbf X_{2n-1}')}}),
\]
$n=1,2,\ldots N$.\footnote{Note that, thanks to the conditioning on the $({\mathbf X_{n}},R_n)$s,
the $J_{n}({\mathbf x})$s may be chosen in such a way that 
$R'_{2n}=0$ for all $n$, but also in such a way that $J_{n}({\mathbf x})>I_N$ for all $n=1,2,\ldots,N$ and all ${\mathbf x}$, so that $R'_{2n}=1$ too may occur. While in the former case no diseased individuals are found among the controls, in the latter the number of diseased and exposed individuals with background characteristics ${\mathbf x}$ among the controls may very well follow the probabilities
${\mathbf P}(T=t,R=r|{\mathbf X}={\mathbf x})$. Whatever the case, if not the whole sample of controls then at least a subset of it can be used to estimate $\mathbf{P}(T=t|R=0,\mathbf{X}=\mathbf{x})$, even if the sample itself does not follow these probabilities (it does in the first case but not in the second).}

The set $({\mathbf X_{1}'},T'_{1},R'_{1}),({\mathbf X_{2}'},T'_{2},R'_{2}),\dots,
({\mathbf X_{2N}'},T'_{2N},R'_{2N})$ so defined represents the data that can be observed in the study. It is no random sample since in it a case is always followed by its control, but more important is the fact that its elements, being only incompletely and selectively observed, do not have the same distribution as the $({\mathbf X_n},T_n,R_n)$s, so it is not obvious that a case-control study can help us investigate the effect of the exposure on the disease. 
To show that it can, note first that, thanks to the way in which the sets 
$\{I_1,I_2,\ldots,I_N\}$ and
$\{J_{1}({\mathbf x}),J_{2}({\mathbf x}),\ldots\}$ are generated,
the data from the cases and the data from the controls, respectively, allow us to estimate the probabilities
\begin{equation}\label{EquationProbabilitiesWeCanEstimate}
p_\mathbf{x} :=\mathbf{P}(T=1|R=1,\mathbf{X}=\mathbf{x} )\quad\mbox{and}\quad q_\mathbf{x} :=\mathbf{P}(T=1|R=0,\mathbf{X}=\mathbf{x} ),
\end{equation}
at least for the $\mathbf{x}\in{\cal X}$ assumed by the observed
${\mathbf X_{n}'}$s.

Of course, these probabilities are not our object of interest: we know that under our model the effect of the exposure on the disease is fully characterized by the probabilities
$\mathbf{P}(R=1|T=t,\mathbf{X}=\mathbf{x})$ as functions of $t$, that is by 
\begin{equation}\label{EquationProbabilitiesOfInterest}
\mathbf{P}(R=1|T=1,\mathbf{X}=\mathbf{x})\quad\mbox{and}\quad \mathbf{P}(R=1|T=0,\mathbf{X}=\mathbf{x}),
\end{equation}
or by some measure of discrepancy between them.
But these probabilities can hardly
be estimated from the data collected in a case-control study because in the set of cases the number
of diseased individuals is fixed and the set of controls, which may have had its source of diseased individuals even more depleted by the preceding sampling of cases, will rarely have a diseased individual to show.

Surprisingly, there is a measure of discrepancy between the probabilities in (\ref{EquationProbabilitiesOfInterest}) which can be expressed in terms of those in
(\ref{EquationProbabilitiesWeCanEstimate}), and hence estimated from the data collected in a case-control study: the \textit{odds ratio}, or, more precisely, the odds ratio conditional on $\{\mathbf{X}=\mathbf{x}\}$, defined by
\begin{equation}\label{Equation4dot5}
{\rm e}(\mathbf{x})=\frac{{\cal O}(R=1|T=1,\mathbf{X}=\mathbf{x})}{{\cal O}(R=1|T=0,\mathbf{X}=\mathbf{x})} ,
\end{equation}
where
\[
{\cal O}(R=1|T=t,\mathbf{X}=\mathbf{x})=\frac{\mathbf{P}(R=1|T=t,\mathbf{X}=\mathbf{x})}{1-\mathbf{P}(R=1|T=t,\mathbf{X}=\mathbf{x})},
\]
the \textit{odds} of an individual with characteristics $\mathbf{x}$
being diseased given that he was exposed to treatment $t$, is an
increasing function of $\mathbf{P}(R=1|T=t,\mathbf{X}=\mathbf{x})$. If the exposure does not cause the disease,
${\cal O}(R=1|T=t,\mathbf{X}=\mathbf{x})$ is constant in
 $t$ and hence ${\rm e}(\mathbf{x})\!=\!1$; otherwise,
$(R=1|T=1,\mathbf{X}=\mathbf{x} )\!>\!(R=1|T=0,\mathbf{X}=\mathbf{x})$
and hence ${\rm e}(\mathbf{x} )>1$. Since, as is easy to verify,
\begin{equation}\label{Equation4dot6}
{\rm e}(\mathbf{x})=\frac{{\cal O}(T=1|R=1,\mathbf{X}=\mathbf{x})}{{\cal O}(T=1|R=0,\mathbf{X}=\mathbf{x})}=
\frac{p_\mathbf{x}(1-q_\mathbf{x})}{q_\mathbf{x}(1-p_\mathbf{x})},
\end{equation}
${\rm e}(\mathbf{x})$ also compares the odds of an individual
having been exposed given that he is diseased with the odds of an
individual having being exposed given that he is not diseased. Thus, although
it is (\ref{Equation4dot5}) that characterizes the causal effect of the exposure on the disease, (\ref{Equation4dot6}) shows that the effect can be estimated from estimates of the 
probabilities (\ref{EquationProbabilitiesWeCanEstimate}).\footnote{The significance of the odds
ratio in case-control studies was first noticed by Cornfield \cite{Cornfield1951}. There exist various types of case-control studies; in particular, several controls may be matched to each case (see, for example, pp.~\!7, 83--86 of \cite{Rosenbaum2002}, where case-control studies are called case-referent studies).}

If we integrate ${\rm e}(\mathbf{x})$ with respect to $\mathbf{P}(\mathbf{X}\le \mathbf{x}|R=1)$ we get an overall measure of the effect of exposure on disease which in principle can be estimated from the observed data, namely
$\mathbf{E}[{\rm e}(\mathbf{X})|R=1]$,
where the conditioning on $\{R=1\}$ reflects the fact that, by virtue of
the sampling scheme, all the values of the $\mathbf{X}_n$s that actually turn up arise conditionally on that event.

\subsection{Determination of a causal effect with an instrumental variable}
\label{SectionInstrumentalVariables}

In situations where unconfoundedness does not hold in the basic model (\ref{CausalModelII}) the study of a causal effect in the full sense of the word is usually ruled out, but sometimes it is still possible to estimate a sort of causal effect. Angrist, Imbens and Rubin \cite{AngristEtAl1996} consider a special situation where the use of an `instrumental variable' permits the estimation of the causal effect in a subpopulation of the
population of interest.\footnote{Some real-life problems in which the method of instrumental variables is potentially useful are described in \cite{ShehanEtAl2008}.}

The model studied in \cite{AngristEtAl1996} is
\begin{equation}\label{ExampleIV1}
R=\rho(\varepsilon,T),\quad T=\tau(\delta,I),
\end{equation}
where the functions $\rho$ and $\tau$ are regarded as unknown, $T$ and $I$ are observable binary random variables, $(\delta,\varepsilon)$ is an unobservable random vector with arbitrary distribution, $I$ and $(\delta,\varepsilon)$ are independent but $\delta$ and $\varepsilon$ are dependent. The model arises in the study of the effect of a treatment, $T$,
on a response, $R$, with the help of a third variable, $I$, called an \emph{instrument}, which is not of direct interest and affects $R$ only through $T$.

In the real-life problem treated in \cite{AngristEtAl1996} the question of interest is whether serving in the military in times of war has detrimental effects on the health of individuals; $R$ stands for the health outcome of a generic individual, $I$ for the individual's
`draft status'---whether he was called to serve or not---and $T$ for his indicator of military service---whether he actually served or not. During the periods covered by the data, recruitment had been determined by a lottery, `low' lottery numbers leading to drafting. Individuals who were called to serve may or may not have avoided joining the army, so the data on such an individual might be $(I,T)=(1,0)$ or $(I,T)=(1,1)$; and individuals who were not called may have volunteered and joined the military, so $(I,T)=(0,1)$ too might be observed. Although in principle there might be individuals who would have refused to serve if recruited but accepted to serve if not recruited, their $\delta$ satisfying $1=\tau(\delta,0)>\tau(\delta,1)=0$, such individuals should be rare; accordingly, it is assumed in \cite{AngristEtAl1996}, and will be assumed here as well, that \emph{$\tau$ is monotonic in its second argument}, i.e.~\!that $\tau(\delta,0)\le \tau(\delta,1)$ for all $\delta$.

In order to study the effect of the treatment on the response---to study how varying
$t$ changes the law of $\rho(\varepsilon,t)$---one might think of estimating
\[
{\cal L}(R|T=t)={\cal L}(\rho(\varepsilon,t)|\tau(\delta,I)=t)
\]
or
\[
{\cal L}(R|T=t,I=i)={\cal L}(\rho(\varepsilon,t)|\tau(\delta,i)=t,I=i),
\]
or functionals of them, from data on $(I,T,R)$. However, the dependence between $\delta$ and $\varepsilon$ shows that the conditioning on $\{\tau(\delta,I)=t\}$ or on $\{\tau(\delta,i)=t,I=i\}$ implies that changes in $t$ cause changes in $\rho(\varepsilon,t)$ also through the first argument of $\rho$, so neither of these probability laws (and indeed no other probability law) describes the effect of the treatment on the response.

Surprisingly, the presence of the `instrument' $I$ in the second equation of (\ref{ExampleIV1}) makes it possible to describe the effect of $T$ on $R$ in some, albeit incomplete, sense.\footnote{The role of $I$ is quite different from that of the vector of potential confounders in model (\ref{CausalModelII}). In fact, one may include such a vector ${\mathbf{X}}$ in both equations of (\ref{ExampleIV1}) and carry out a stratified version of the analysis that follows by conditioning on $\{{\mathbf{X}}={\mathbf{x}}\}$.}
Indeed, Angrist, Imbens and Rubin show that the parameter
\begin{equation}\label{ExampleIV1A}
\theta:=\frac{{\mathbf{E}}(R|I=1)-{\mathbf{E}}(R|I=0)}
{{\mathbf{E}}(T|I=1)-{\mathbf{E}}(T|I=0)},
\end{equation}
which is defined whenever $I$ and $T$ are correlated, and which can be estimated from data on $(I,T,R)$ by replacing expected values by sample averages, is a measure of the causal effect of $T$ on $R$.

To see this, note that by the independence of $I$ and $\delta$ and the monotonicity of $\tau$,
which implies
\[
\tau(\delta,1)-\tau(\delta,0)\neq 0\,\,\,\,\Leftrightarrow\,\,\,\,
\tau(\delta,1)-\tau(\delta,0)=1\,\,\,\,\Leftrightarrow\,\,\,\, 
\tau(\delta,1)=1\,\,\wedge\,\,\tau(\delta,0)=0,
\]
the denominator of $\theta$ is
\begin{eqnarray*}
{\mathbf{E}}(T|I=1)-{\mathbf{E}}(T|I=0)\!\!\!&\!=\!&\!\!\!
{\mathbf{E}}[\tau(\delta,1)|I=1]-{\mathbf{E}}[\tau(\delta,0)|I=0]\\
\!\!\!&\!=\!&\!\!\!
{\mathbf{E}}[\tau(\delta,1)-\tau(\delta,0)]\\
\!\!\!&\!=\!&\!\!\!
{\mathbf{P}}[\tau(\delta,1)-\tau(\delta,0)=1].
\end{eqnarray*}
By the independence of $I$ and $(\delta,\varepsilon)$ and by  
\[
\rho\left(\varepsilon,\tau(\delta,1)\right)-\rho\left(\varepsilon,\tau(\delta,0)\right)\neq 0
\,\,\,\,\Rightarrow\,\,\,\,
\tau(\delta,1)>\tau(\delta,0),
\]
the numerator is
\begin{eqnarray*}
{\mathbf{E}}(R|I=1)-{\mathbf{E}}(R|I=0)\!\!\!&\!=\!&\!\!\!
{\mathbf{E}}[\rho\left(\varepsilon,\tau(\delta,1)\right)|I=1]-
{\mathbf{E}}[\rho\left(\varepsilon,\tau(\delta,0)\right)|I=0]\\
\!\!\!&\!=\!&\!\!\!
{\mathbf{E}}[\rho\left(\varepsilon,\tau(\delta,1)\right)-\rho\left(\varepsilon,\tau(\delta,0)\right)]\\
\!\!\!&\!=\!&\!\!\!
{\mathbf{E}}[\left\{\rho\left(\varepsilon,1\right)-\rho\left(\varepsilon,0\right)\right\}{\mathbf{1}}_{\{\tau(\delta,1)-\tau(\delta,0)=1\}}].
\end{eqnarray*}
Thus
\[
\theta=
\frac{{\mathbf{E}}[\left\{\rho\left(\varepsilon,1\right)-\rho\left(\varepsilon,0\right)\right\}{\mathbf{1}}_{\{\tau(\delta,1)-\tau(\delta,0)=1\}}]}
{{\mathbf{P}}[\tau(\delta,1)-\tau(\delta,0)=1]}
\]
\vspace{-0.6cm}
\begin{equation}\label{ExampleIV1B}
\mbox{\quad}
\end{equation}
\vspace{-0.6cm}
\[
\,\,\,\,\,\,\,\,={\mathbf{E}}\left[\rho\left(\varepsilon,1\right)-\rho\left(\varepsilon,0\right)
\left|\tau(\delta,1)-\tau(\delta,0)=1\right.\right]
\]
is the \emph{average effect of treatment on response in the subpopulation of individuals for whom $\tau(\delta,1)-\tau(\delta,0)=1$}.\footnote{In \cite{AngristEtAl1996}, where $I$ is really a suggestion of treatment rather than the treatment $T=\tau(\delta,I)$, such individuals are called `compliers': they get treatment only if treatment is suggested to them.}

Although it is impossible to identify the subpopulation of individuals to which the causal effect represented by $\theta$ applies ($\tau(\delta,1)-\tau(\delta,0)$ being unobservable), evidence that $\theta\neq 0$ does provide evidence of a treatment effect.

\begin{remarks} {\bf (i)} The study of the \emph{existence} of the treatment effect is tackled very simply by writing (\ref{ExampleIV1B}) in the form of (\ref{CausalModelII}) with unconfoundedness, namely as
\[
R=\tilde\rho(V,I)\equiv \rho(\varepsilon,\tau(\delta,I)),\quad I=\iota(U),
\]
where $\iota$ is some function and $U$ and $V$ are independent standard uniforms, the latter being obtained by alternating the digits in the decimal expansions of $\delta$ and $\varepsilon$ (assumed uniform without loss of generality). Regarding the instrument as a treatment, the treatment effect is characterized by
${\cal L}(R|I=i)={\cal L}(\tilde\rho(V,i))$ for varying $i$, and if evidence is found of an effect of $I$ on $R$ then it follows from the assumed correlation between $I$ and $T=\tau(\delta,I)$ and the assumption that $T$ exerts its influence on $R$ through the second argument of $\rho$ that the evidence is equally applicable to the effect of $T$ on $R$.

This is true without any assumptions about the range of $R$ and $T$; moreover,
\[
\mathbf{E}\left[R\left|I=i\right.\right]-
\mathbf{E}\left[R\left|I=i'\right.\right]=
\mathbf{E}\left[\rho(\varepsilon,\tau(\delta,i))-\rho(\varepsilon,\tau(\delta,i'))\right]=
\]
\[
\mathbf{E}\left[\rho(\varepsilon,\tau(\delta,i))-\rho(\varepsilon,\tau(\delta,i'))\left|
\tau(\delta,i)\neq \tau(\delta,i')\right.\right]\mathbf{P}\left(\tau(\delta,i)\neq \tau(\delta,i')\right)\!,
\]
so what one may consider to be the effect of the treatment of interest, albeit in an unobservable  subpopulation characterized by the condition $\tau(\delta,i)\neq \tau(\delta,i')$, is proportional to the effect of the instrument on the response. What is gained by making additional assumptions (such as $\tau$ being monotonic and $R$ and $T$ being binary) is the possibility of identifying the proportionality constant as something that can be estimated from observed data, in the case above as $\mathbf{P}\left(\tau(\delta,1)\neq \tau(\delta,0)\right)={\mathbf{E}}(T|I=1)-{\mathbf{E}}(T|I=0)$.

{\bf (ii)} If $T$ is discrete rather than binary, similar arguments show that
\begin{eqnarray*}
\theta\!\!\!&\!=\!&\!\!\!
\frac{{\mathbf{E}}[\left\{\rho\left(\varepsilon,\tau(\delta,1)\right)-\rho\left(\varepsilon,\tau(\delta,0)\right)\right\}{\mathbf{1}}_{\{\tau(\delta,1)>\tau(\delta,0)\}}]}
{{\mathbf{E}}[\left\{\tau(\delta,1)-\tau(\delta,0)\right\}{\mathbf{1}}_{\{\tau(\delta,1)>\tau(\delta,0)\}}]}\\
\!\!\!&\!\mbox{\quad}\!&\!\!\!\vspace{-1cm}\\
\!\!\!&\!=\!&\!\!\!
\frac{{\mathbf{E}}\left[\rho\left(\varepsilon,\tau(\delta,1)\right)-\rho\left(\varepsilon,\tau(\delta,0)\right)
\left|\tau(\delta,1)>\tau(\delta,0)\right.\right]}
{{\mathbf{E}}\left[\tau(\delta,1)-\tau(\delta,0)
\left|\tau(\delta,1)>\tau(\delta,0)\right.\right]}.
\end{eqnarray*}
This, too, is a measure of the effect of the treatment on the response in the subpopulation of individuals for whom $\tau(\delta,1)>\tau(\delta,0)$; it is more difficult to interpret than the first because one cannot estimate the last denominator, but, again,
evidence that $\theta\neq 0$ provides evidence of a treatment effect.

{\bf (iii)}
The analysis can be extended to model (\ref{ExampleIV1}) in the case where $I$ and $T$ are discrete random variables with $I$ taking the numerical values $i_0<i_1<i_2<\cdots$ and $\tau$ such that $\tau(\delta,i_0)\le \tau(\delta,i_k)$ for $k\ge 1$. For simplicity, assume that
${\mathbf{P}}[\tau(\delta,i_k)\!>\!\tau(\delta,i_0)]\!>\!0$ for all $k$ and consider the parameters
\[
\theta_k:=\frac{{\mathbf{E}}(R|I=i_k)-{\mathbf{E}}(R|I=i_0)}
{{\mathbf{E}}(T|I=i_k)-{\mathbf{E}}(T|I=i_0)}=
\frac{{\mathbf{E}}[\rho\left(\varepsilon,\tau(\delta,i_k)\right)-\rho\left(\varepsilon,\tau(\delta,i_0)\right)]}{{\mathbf{E}}[\tau(\delta,i_k)-\tau(\delta,i_0)]},
\]
\[
\,=
\frac{
{\mathbf{E}}\left[\rho\left(\varepsilon,\tau(\delta,i_k)\right)-\rho\left(\varepsilon,\tau(\delta,i_0)\right)
\left|\tau(\delta,i_k)>\tau(\delta,i_0)\right.\right]
}
{
{\mathbf{E}}\left[\tau(\delta,i_k)-\tau(\delta,i_0)
\left|\tau(\delta,i_k)>\tau(\delta,i_0)\right.\right]},
\quad\quad\quad\quad
\]
which can be estimated from data on $(I,T,R)$ and measure the effect of the treatment in the subpopulations of individuals for whom $\tau(\delta,i_k)>\tau(\delta,i_0)$.
The numbers
\[
p_k :=\frac{{\mathbf{P}}(I=i_k)\left\{{\mathbf{E}}(T|I=i_k)-{\mathbf{E}}(T|I=i_0)\right\}}
{\sum\nolimits_{l\ge 1}
{\mathbf{P}}(I=i_l)\left\{{\mathbf{E}}(T|I=i_l)-{\mathbf{E}}(T|I=i_0)\right\}
}
\]
define a probability distribution on $\mathbb{N}$ which can be estimated from data on
$(I,T)$, so we can define an overall measure of treatment effect by
$\Theta:=\sum\nolimits_{k\ge 1}\theta_k p_k$, which turns out to be
\begin{eqnarray*}
\Theta
=
\frac{
{\mathbf{E}}[\rho\left(\varepsilon,\tau(\delta,I)\right)-\rho\left(\varepsilon,\tau(\delta,i_0)\right)\left|\right. \tau(\delta,I)>\tau(\delta,i_0)]}
{{\mathbf{E}}[\tau(\delta,I)-\tau(\delta,i_0)\left|\right. \tau(\delta,I)>\tau(\delta,i_0)]}.
\end{eqnarray*}
Thus $\Theta$, which can be estimated from data,
is a measure (though not a very clear one) of the effect of treatment on response in the subpopulation of individuals for whom $\tau(\delta,I)>\tau(\delta,i_0)$. When $T$ is binary it simplifies to a conditional expectation and compares with (\ref{ExampleIV1B}).\quad\quad\quad\quad\quad\quad\quad\quad\quad\quad\quad\quad\quad\quad\quad\quad\quad$\square$
\end{remarks}

\section{Pearl's calculus of intervention}\label{CalculusOfIntervention}

As suggested by N. Singpurwalla (see \cite{Singpurwalla2002} and the ensuing exchange with Pearl), for example, Pearl's calculus of intervention ought to admit a formulation based entirely on notation and results of elementary probability theory (of the kind used throughout our work). We believe this to be the case,  and in this section provide an interpretation and proofs of the first two of Pearl's calculus rules; the third rule is beyond our grasp.

We speak of an {\it interpretation} because we are not
sure that the rules stated here correspond to Pearl's; see pp.~\!85--86 of \cite{Pearl2009a}
for the statements of the rules and \cite{Pearl1995a} for proofs---the only ones we have seen and which we are unable to understand. As far as we see---and if we exclude the efforts by R.~Tucci in \cite{Tucci2013}, which we are equally unable to understand---no additional elucidation of the calculus and no alternative proofs of its rules have been given by other authors, nor by Pearl (not even in \cite{Pearl2012}), since their publication in \cite{Pearl1995a}; but the three rules have been widely cited and reproduced almost verbatim from this last source.

Consider a second-level model ${\cal M}$ consisting of four disjoint sets of
variables---also called `nodes' in connection with the graph representing the relationships between them---${\mathbf W}$, ${\mathbf X}$, ${\mathbf Y}$ and ${\mathbf Z}$.
Regarding ${\mathbf W}$, ${\mathbf X}$, ${\mathbf Y}$ and ${\mathbf Z}$ as random vectors, such a model is defined by a system of equations of the form
\begin{equation}\label{ModelM}
\left\{
\begin{array}{l}
{\mathbf W}=\varphi_1(U_1,{\mathbf W},{\mathbf X},{\mathbf Y},{\mathbf Z}), \\
{\mathbf X}=\varphi_2(U_2,{\mathbf W},{\mathbf X},{\mathbf Y},{\mathbf Z}), \\
{\mathbf Y}=\varphi_3(U_3,{\mathbf W},{\mathbf X},{\mathbf Y},{\mathbf Z}), \\
{\mathbf Z}=\varphi_4(U_4,{\mathbf W},{\mathbf X},{\mathbf Y},{\mathbf Z}),
\end{array}
\right.
\end{equation}
where the $U_i$s are independent vectors of independent uniforms and the $\varphi_i$s are vector-valued functions. In such a system, a variable appearing on the left-hand
side of an equality is understood not to enter as a variable in the corresponding coordinate of the function on the right. For instance, if $W_j$ is the $j$-th coordinate of ${\mathbf W}$ then $W_j$ plays no active role in the $j$-th coordinate function of $\varphi_1$.
More generally, a variable defined as a function of another variable cannot in turn be involved in its definition.
Thus the first line in (\ref{ModelM}) could be something like
\begin{equation*}
\left\{
\begin{array}{l}
W_1=-\log{U_{1,1}},\\
W_2=1/(1+U_{1,2}),\\
W_3=a_1\exp\left({a_2 U_{2,1}+a_3 W_2 X_2 + a_4 W_1 Y_6 Z_2}\right),\\
W_4=b_1\sin\left({U_{2,2}W_1 W_3}\right)+b_2\log\left({Z_1^2+b_3 W_2  e^{Y_3}}\right),\\
\end{array}
\right.
\end{equation*}
with $U_{1,1}$, $U_{1,2}$ coordinates of $U_{1}$,
$U_{2,1}$, $U_{2,2}$ coordinates of $U_{2}$, $W_1,W_2,W_3,W_4$ coordinates of ${\mathbf W}$,
$X_2$ of ${\mathbf X}$, etc.; and in this case $X_2$ could not be a function of $(W_3,W_4)$ 
because $W_3$ is a function of $X_2$ and $W_4$ a function of $W_3$.

In order to formulate Pearl's rules we need to distinguish between 
exogenous variables---i.e.~those that are functions of uniforms alone, such as $W_1$ and $W_2$ in the example just given---from the non-exogenous ones, so we write (\ref{ModelM}) as
\begin{equation}\label{ModelM2}
\left\{
\begin{array}{l}
\!{\mathbf W}_1=\varphi_{1}(U_1),\quad
\!{\mathbf W}_2=\phi_{1}(V_1,{\mathbf W},{\mathbf X},{\mathbf Y},{\mathbf Z}), \\
{\mathbf X}_1=\varphi_{2}(U_2),\quad
{\mathbf X}_2=\phi_{2}(V_2,{\mathbf W},{\mathbf X},{\mathbf Y},{\mathbf Z}), \\
{\mathbf Y}_1=\varphi_{3}(U_3),\quad
{\mathbf Y}_2=\phi_{3}(V_3,{\mathbf W},{\mathbf X},{\mathbf Y},{\mathbf Z}), \\
\,{\mathbf Z}_1=\varphi_{4}(U_4),\quad\,
{\mathbf Z}_2=\phi_{4}(V_4,{\mathbf W},{\mathbf X},{\mathbf Y},{\mathbf Z}),
\end{array}
\right.
\end{equation}
with the $U_i$s and $V_i$s independent vectors of independent uniforms, the
$\varphi_i$s and $\phi_i$s vector-valued functions,
${\mathbf W}=({\mathbf W}_1,{\mathbf W}_2)$,
${\mathbf X}=({\mathbf X}_1,{\mathbf X}_2)$, etc.

Pearl's rules are identities, valid under conditions to be stated below, between certain conditional probability functions pertaining to two intervention models derived from ${\cal M}$: a model ${\cal M}_{\mathbf{x}}'$ corresponding to a numerical vector
${\mathbf x}=({\mathbf x}_1,{\mathbf x}_2)$ in the range of ${\mathbf X}$ and defined by
\begin{equation}\label{ModelM2prime}
\left\{
\begin{array}{l}
\!{\mathbf W}'_{1}=\varphi_{1}(U_1),\quad
\!{\mathbf W}'_{2,{\mathbf x}}=\phi_{1}(V_1,{\mathbf W}'_{{\mathbf x}},{\mathbf x}_1,{\mathbf x}_2,{\mathbf Y}'_{{\mathbf x}},{\mathbf Z}'_{{\mathbf x}}), \\
{\mathbf Y}'_{1}=\varphi_{3}(U_3),\quad
{\mathbf Y}'_{2,{\mathbf x}}=\phi_{3}(V_3,{\mathbf W}'_{{\mathbf x}},{\mathbf x}_1,{\mathbf x}_2,{\mathbf Y}'_{{\mathbf x}},{\mathbf Z}'_{{\mathbf x}}), \\
\,{\mathbf Z}'_{1}=\varphi_{4}(U_4),\quad\,
{\mathbf Z}'_{2,{\mathbf x}}=\phi_{4}(V_4,{\mathbf W}'_{{\mathbf x}},{\mathbf x}_1,{\mathbf x}_2,{\mathbf Y}'_{{\mathbf x}},{\mathbf Z}'_{{\mathbf x}}),
\end{array}
\right.
\end{equation}
where we write
${\mathbf W}'_{{\mathbf x}}=({\mathbf W}'_{1},{\mathbf W}'_{2,{\mathbf x}})$,
${\mathbf Y}'_{{\mathbf x}}=({\mathbf Y}'_{1},{\mathbf Y}'_{2,{\mathbf x}})$ and
${\mathbf Z}'_{{\mathbf x}}=({\mathbf Z}'_{1},{\mathbf Z}'_{2,{\mathbf x}})$;
and a model ${\cal M}_{\mathbf{x},\mathbf{z}}''$ corresponding to numerical vectors
${\mathbf x}=({\mathbf x}_1,{\mathbf x}_2)$ and ${\mathbf z}$ in the ranges of
${\mathbf X}$ and ${\mathbf Z}$, and defined by
\begin{equation}\label{ModelM2primeprime}
\left\{
\begin{array}{l}
\!{\mathbf W}''_{1}=\varphi_{1}(U_1),\quad
\!{\mathbf W}''_{2,{\mathbf x},{\mathbf z}}=\phi_{1}(V_1,{\mathbf W}''_{{\mathbf x},{\mathbf z}},{\mathbf x}_1,{\mathbf x}_2,{\mathbf Y}''_{{\mathbf x},{\mathbf z}},{\mathbf z}), \\
{\mathbf Y}''_{1}=\varphi_{3}(U_3),\quad
{\mathbf Y}''_{2,{\mathbf x},{\mathbf z}}=\phi_{3}(V_3,{\mathbf W}''_{{\mathbf x},{\mathbf z}},{\mathbf x}_1,{\mathbf x}_2,{\mathbf Y}''_{{\mathbf x},{\mathbf z}},{\mathbf z}), 
\end{array}
\right.
\end{equation}
where ${\mathbf W}''_{{\mathbf x},{\mathbf z}}=({\mathbf W}''_{1},{\mathbf W}''_{2,{\mathbf x},{\mathbf z}})$ and
${\mathbf Y}''_{{\mathbf x},{\mathbf z}}=({\mathbf Y}''_{1},{\mathbf Y}''_{2,{\mathbf x},{\mathbf z}})$.

Note that ${\cal M}_{\mathbf{x}}'$ is obtained from ${\cal M}$ by removing all the equations that define nodes of ${\mathbf {X}}$ 
and replacing the nodes of ${\mathbf X}$ in all the remaining equations by elements of a numerical vector $\mathbf{x}$, and
${\cal M}_{\mathbf{x},\mathbf{z}}''$ is obtained from ${\cal M}$ by removing all the equations that define nodes of ${\mathbf {X}}$ and ${\mathbf Z}$ 
and replacing the nodes of ${\mathbf X}$ and ${\mathbf Z}$ in all the remaining equations by elements of numerical vectors $\mathbf{x}$ and $\mathbf{z}$.

\subsection{Rule 1}\label{AppendixRule1}

The first identity is
\[
P(y|{\mathbf{\hat x}},z,w):=
{\mathbf P}({\mathbf Y}_{\mathbf{x}}'=y|{\mathbf Z}_{\mathbf{x}}'=z,{\mathbf {W}}_{\mathbf{x}}'=w)
\quad\quad\quad\quad\quad\quad\quad\quad\quad
\]
\vspace{-0.8cm}
\begin{equation}\label{Rule1}
\mbox{\quad}
\end{equation}
\vspace{-0.8cm}
\[
\quad\,
={\mathbf P}({\mathbf Y}_{\mathbf{x}}'=y|{\mathbf {W}}_{\mathbf{x}}'=w)
=:P(y|{\mathbf{\hat x}},w),
\quad\quad
\]
%
where the right- and leftmost terms indicate, in Pearl's notation,
conditional probabilities pertaining to the intervention model ${\cal M}_{\mathbf{x}}'$
 (the identity proper is the middle equality) and
$w,y,z$, like the $\mathbf{x}$ and $\mathbf{z}$ that determine ${\cal M}_{\mathbf{x}}'$ and ${\cal M}_{\mathbf{x},\mathbf{z}}''$, denote numerical vectors.\footnote{In Pearl's notation a symbol such as ${\mathbf{\hat x}}$ serves to indicate the model from which the probabilities are to be computed (in this case
${\cal M}_{\mathbf{x}}'$, which is determined by $\mathbf{x}$);
sometimes ${\mathbf{\hat x}}$ is
replaced by $do({\mathbf{x}})$, where the $do$ stands for the operation that 
transforms a second-level model into an intervention version of it.}

It is valid under the following condition: 

\vspace{0.25cm}
{\bf C1.}~Let the equation ${\mathbf{X}}'_1=\varphi_{2}(U_2)$ be added to
${\cal M}_{\mathbf{x}}'$; then
\[
{\mathbf P}\left({\mathbf Y}_{\mathbf{x}}'=y\left|
{\mathbf {W}}_{\mathbf{x}}'=w,
{\mathbf{X}}'_1={\mathbf{x}}_1,
{\mathbf {Z}}_{\mathbf{x}}'=z
\right.\right)=
{\mathbf P}\left({\mathbf Y}_{\mathbf{x}}'=y\left|
{\mathbf {W}}_{\mathbf{x}}'=w,
{\mathbf{X}}'_1={\mathbf{x}}_1
\right.\right)
\]
for all $(w,y,z)$ in the range of
$({\mathbf {W}}_{\mathbf{x}}',{\mathbf Y}_{\mathbf{x}}',{\mathbf {Z}}_{\mathbf{x}}')$.
\vspace{0.25cm}

Rule 1 follows in three steps from the exogeneity of ${\mathbf{X}}'_1$, from {\bf C1} and again from the exogeneity of ${\mathbf{X}}'_1$:
\[
{\mathbf P}({\mathbf Y}_{\mathbf{x}}'=y|{\mathbf Z}_{\mathbf{x}}'=z,{\mathbf {W}}_{\mathbf{x}}'=w)=
\frac{{\mathbf P}({\mathbf Z}_{\mathbf{x}}'=z,{\mathbf Y}_{\mathbf{x}}'=y,
{\mathbf {W}}_{\mathbf{x}}'=w){\mathbf P}({\mathbf X}'_1={\mathbf x}_1)}
{{\mathbf P}({\mathbf Z}_{\mathbf{x}}'=z,{\mathbf {W}}_{\mathbf{x}}'=w)
{\mathbf P}({\mathbf X}'_1={\mathbf x}_1)}\quad\quad
\]
\[
\quad\quad\quad\quad\quad\quad\quad\quad\,\,\,\,\,\,\,\,\,\,\!
=\frac{{\mathbf P}({\mathbf Z}_{\mathbf{x}}'=z,{\mathbf Y}_{\mathbf{x}}'=y,
{\mathbf X}'_1={\mathbf x}_1,{\mathbf {W}}_{\mathbf{x}}'=w)}
{{\mathbf P}({\mathbf Z}_{\mathbf{x}}'=z,{\mathbf X}_1={\mathbf x}_1,
{\mathbf {W}}_{\mathbf{x}}'=w)}
\]
\[
\quad\quad\quad\quad\quad\quad\quad\,\,\,\,\,\,\,\,\,\,\,\,\,\!
={\mathbf P}({\mathbf Y}_{\mathbf{x}}'=y|{\mathbf Z}_{\mathbf{x}}'=z,
{\mathbf X}'_1={\mathbf x}_1,{\mathbf {W}}_{\mathbf{x}}'=w)
\]
\[
\quad\quad\quad\quad\quad\,\,\,\,\,\!\!
={\mathbf P}({\mathbf Y}_{\mathbf{x}}'=y|
{\mathbf X}'_1={\mathbf x}_1,{\mathbf {W}}_{\mathbf{x}}'=w)
\]
\[
\quad\quad\,\!\!\!
={\mathbf P}({\mathbf Y}_{\mathbf{x}}'=y|{\mathbf {W}}_{\mathbf{x}}'=w).
\]

\subsection{Rule 2}\label{AppendixRule2}

The second identity is
\[
P(y|{\mathbf {\hat x}},{\mathbf {\hat z}},w):=
{\mathbf P}({\mathbf Y}_{\mathbf{x},\mathbf{z}}''=y|{\mathbf W}_{\mathbf{x},\mathbf{z}}''=w)
\quad\quad\quad\quad\quad\quad\quad\quad\quad
\]
\vspace{-0.8cm}
\begin{equation}\label{Rule2}
\mbox{\quad}
\end{equation}
\vspace{-0.8cm}
\[
\quad\quad\quad\quad\quad\quad\quad\quad\!
={\mathbf P}({\mathbf Y}_{\mathbf{x}}'=y|{\mathbf Z}_{\mathbf{x}}'={\mathbf z},
{\mathbf {W}}_{\mathbf{x}}'=w)
=:P(y|{\mathbf {\hat x}},\mathbf {z},w),
\quad\quad
\]
in notation similar to that of (\ref{Rule1}).

It is valid under the following condition: 

\vspace{0.25cm}
{\bf C2.}~Let the equations ${\mathbf{X}}''_1=\varphi_{2}(U_2)$,
${\mathbf{Z}}''_{1}=\varphi_{4}(U_4)$ and
\[
{\mathbf{Z}}''_{2,{\mathbf x}}=
\phi_{4}(V_4,{\mathbf{W}}''_{{\mathbf x},{\mathbf z}},{\mathbf x}_1,{\mathbf x}_2,
{\mathbf{Y}}''_{{\mathbf x},{\mathbf z}},{\mathbf{Z}}''_{{\mathbf x}})
\]
be added to ${\cal M}_{{\mathbf{x}},{\mathbf{z}}}''$, and write
${\mathbf{Z}}''_{{\mathbf x}}=
({\mathbf{Z}}''_{1},{\mathbf{Z}}''_{2,{\mathbf x}})$; then
\[
{\mathbf P}\left({\mathbf{Y}}''_{{\mathbf x},{\mathbf z}}\!=\!y\left|
{\mathbf{W}}''_{{\mathbf x},{\mathbf z}}\!=\!w,
{\mathbf{X}}''_1\!=\!{\mathbf{x}}_1,
{\mathbf{Z}}''_{{\mathbf x}}\!=\!{\mathbf z}
\right.\right)\!=\!
{\mathbf P}\left({\mathbf{Y}}''_{{\mathbf x},{\mathbf z}}\!=\!y\left|
{\mathbf{W}}''_{{\mathbf x},{\mathbf z}}\!=\!w,
{\mathbf{X}}''_1\!=\!{\mathbf{x}}_1
\right.\right)
\]
for all $(w,y)$ in the range of
$({\mathbf{W}}''_{{\mathbf x},{\mathbf z}},{\mathbf{Y}}''_{{\mathbf x},{\mathbf z}})$.
\vspace{0.25cm}

Note that by (\ref{ModelM2prime}), (\ref{ModelM2primeprime}), and the 
definitions in {\bf C1} and {\bf C2} we have
\[
{\mathbf P}\left({\mathbf{Y}}''_{{\mathbf x},{\mathbf z}}\!=\!y,
{\mathbf{W}}''_{{\mathbf x},{\mathbf z}}\!=\!w,
{\mathbf{X}}''_1\!=\!{\mathbf{x}}_1,
{\mathbf{Z}}''_{{\mathbf x}}\!=\!{\mathbf z}
\right)=
{\mathbf P}\left({\mathbf{Y}}'_{{\mathbf x}}\!=\!y,
{\mathbf{W}}'_{{\mathbf x}}\!=\!w,
{\mathbf{X}}'_1\!=\!{\mathbf{x}}_1,
{\mathbf{Z}}'_{{\mathbf x}}\!=\!{\mathbf z}
\right).
\]
Rule 2 follows from this identity, {\bf C2} and the exogeneity of ${\mathbf{X}}'_1$ and
${\mathbf{X}}''_1$:
\[
{\mathbf P}({\mathbf{Y}}''_{{\mathbf x},{\mathbf z}}=y|
{\mathbf{W}}''_{{\mathbf x},{\mathbf z}}=w)=
\frac{{\mathbf P}({\mathbf{Y}}''_{{\mathbf x},{\mathbf z}}=y,
{\mathbf X}''_1={\mathbf x}_1,
{\mathbf{W}}''_{{\mathbf x},{\mathbf z}}=w)}
{{\mathbf P}({\mathbf X}''_1={\mathbf x}_1,
{\mathbf{W}}''_{{\mathbf x},{\mathbf z}}=w)}\quad\quad
\]
\vspace{-0.25cm}
\[
\quad\quad\quad\quad\quad\quad\quad\,\,\,\,\,\,\,\,\!
={\mathbf P}({\mathbf{Y}}''_{{\mathbf x},{\mathbf z}}=y|{\mathbf X}''_1={\mathbf x}_1,
{\mathbf{W}}''_{{\mathbf x},{\mathbf z}}=w)
\]
\[
\quad\quad\quad\quad\quad\quad\quad\quad\quad\,\,\,\,\,\,\,\,\,\,\,\,\,\!
={\mathbf P}\left({\mathbf{Y}}''_{{\mathbf x},{\mathbf z}}\!=\!y\left|
{\mathbf{W}}''_{{\mathbf x},{\mathbf z}}\!=\!w,
{\mathbf{X}}''_1\!=\!{\mathbf{x}}_1,
{\mathbf{Z}}''_{{\mathbf x}}\!=\!{\mathbf z}
\right.\right)
\]
\[
\quad\quad\quad\quad\quad\quad\quad\quad\,\,\,\,\,\,\,\,\,\,\!
={\mathbf P}\left({\mathbf{Y}}'_{{\mathbf x}}\!=\!y\left|
{\mathbf{W}}'_{{\mathbf x}}\!=\!w,
{\mathbf{X}}'_1\!=\!{\mathbf{x}}_1,
{\mathbf{Z}}'_{{\mathbf x}}\!=\!{\mathbf z}
\right.\right)
\]
\[
\quad\quad\quad\quad\quad\,\,\,\,\,\,\,\!
={\mathbf P}\left({\mathbf{Y}}'_{{\mathbf x}}\!=\!y\left|
{\mathbf{W}}'_{{\mathbf x}}\!=\!w,
{\mathbf{Z}}'_{{\mathbf x}}\!=\!{\mathbf z}
\right.\right)\!.
\]

\end{document}